\documentclass[10pt, a4paper]{article}

\usepackage{amsfonts, latexsym, amsmath, amssymb, fancyhdr, amsthm}
\usepackage[english]{babel}

\newtheorem{theorem}{Theorem}[section]
\newtheorem{lemma}{Lemma}[section]
\newtheorem{proposition}{Proposition}[section]
\newtheorem{definition}{Definition}[section]
\newtheorem{corollary}{Corollary}[section]
\newtheorem{remark}{Remark}

\title{\textbf{Existence, positivity and stability for a nonlinear model of cellular proliferation\thanks{This paper has been published in Nonlinear Analysis: Real World Applications, 6, 337-366, 2005.}}}
\author{Mostafa Adimy\thanks{E-mail: mostafa.adimy@univ-pau.fr} \quad and \quad Fabien Crauste\thanks{E-mail: fabien.crauste@univ-pau.fr}}
\date{Year 2004}

\begin{document}

\maketitle
\begin{center}
\emph{Laboratoire de Math\'ematiques Appliqu\'ees, FRE 2570}\\
\emph{Universit\'e de Pau et des Pays de l'Adour,}\\
\emph{Avenue de l'universit\'e, 64000 Pau, France}
\end{center}

\bigskip{}

\begin{abstract}
In this paper, we investigate a system of two nonlinear partial
differential equations, arising from a model of cellular
proliferation which describes the production of blood cells in the
bone marrow. Due to cellular replication, the two partial
differential equations exhibit a retardation of the maturation
variable and a temporal delay depending on this maturity. We show
that this model has a unique solution which is global under a
classical Lipschitz condition. We also obtain the positivity of
the solutions and the local and global stability of the trivial
equilibrium.
\end{abstract}

\bigskip{}


\section{Introduction}

We analyse, in this paper, a mathematical model arising from the
blood production system. It is based on a system proposed by
Mackey and Rudnicki \cite{mackey1994}, in 1994, to describe the
dynamics of hematopoietic stem cells in the bone marrow. The
origin of this system is a model of Burns and Tannock \cite{burns}
(1970) in which each cell can be either in a proliferating phase
or in a resting phase (also called $G_0$-phase). The resulting
model is a time-age-maturity structured system.

Proliferating cells are in the cell cycle, in that they are
committed to divide at the end of the mitosis, the so-called point
of cytokinesis. After division, they  give birth to two daughter
cells which enter immediatly the resting phase. Proliferating
cells can also die by apoptosis, a programmed cell death.

The resting phase is a quiescent stage in the cellular
development. Cells in this phase can not divide: they mature and,
provided they do not die, they enter the proliferating phase and
complete the cycle.

The model in \cite{mackey1994} has been analysed by Mackey and Rey
\cite{mackey1995_1,mackey1995_2} in 1995, Crabb \emph{et al.}
\cite{crabb1996_1,crabb1996_2} in 1996, Dyson \emph{et al.}
\cite{webb1996} in 1996 and Adimy and Pujo-Menjouet
\cite{adimypujo,adimypujo2} in 2001 and 2003. In these studies,
the authors assumed that all cells divide exactly at the same age.

However, in the most general situation in a cellular population,
it is believed that the time required for a cell to divide is not
identical between cells (see Bradford \emph{et al.}
\cite{bradford}). For example, pluripotent stem cells (which are
the less mature cells) divide faster than committed stem cells,
which are the more mature stem cells. In 1993, Mackey and Rey
\cite{mackey1993} considered a model in which the time required
for a cell to divide is distributed according to a density, but
the authors only made a numerical analysis of their model. Dyson
\emph{et al.} \cite{webb2000, webb2000_2}, in 2000, considered a
time-age-maturity structured equation in which all cells do not
divide at the same age. They presented the basic theory of
existence, uniqueness and properties of the solution operator.
However, in their model, they considered only one phase (the
proliferating one), and the intermediary flux between the two
phases is not represented. In 2003, Adimy and Crauste
\cite{adimycrauste2003} considered a model in which the
proliferating phase duration is distributed according to a density
with compact support. They obtained global stability results for
their model.

In this work, we consider the situation when the age at
cytokinesis depends on the maturity of the cell at the point of
commitment, that means when it enters the proliferating phase. We
assume that each cell entering the proliferating phase with a
maturity $m$ divides at age $\tau=\tau(m)$, depending on this
maturity. This hypothesis can be found, for example, in Mitchison
\cite{mitchison} (1971) and John \cite{john} (1981). This yields
to the boundary condition (\ref{bctauvariant2}). To our knowledge,
nobody has studied this model, except Adimy and Pujo-Menjouet in
\cite{adimypujo2003}, where they considered only a linear case.

We obtain a system of first order partial differential equations
with a time delay depending on the maturity and a retardation of
the maturation variable. We investigate the basic theory of
existence, uniqueness, positivity and stability of the solutions
of our model.

The paper is organised as follows. In Section \ref{presentation},
we present the time-age-maturity structured model. By using the
characteristics method, we reduce this model to a time-maturity
structured system, which is formed by two partial differential
equations with a time delay depending on the maturity and a
nonlocal dependence in the maturity variable. In Section
\ref{localexistence}, we first give an integrated formulation of
our model by using the classical variation of constant formula and
then we prove local existence of solutions, by using a fixed-point
theorem, and their global continuation. We deduce the global
existence. In Section \ref{sectionpositivityregularity}, we obtain
the positivity of these solutions by developping a method
described by Webb \cite{webb1985}. In Section
\ref{sectionlocalstability}, we concentrate on the stability of
the trivial equilibrium of the system and, in the last section, we
discuss the model and the asymptotic behaviour.

\section{Biological background and equations of the model}
\label{presentation}

Each cell is caracterised, in the two phases, by its age and its
maturity. The maturity describes the development of the cell. It
is the concentration of what composes a cell, such as proteins or
other elements one can measure experimentally. The maturity is
supposed to be a continuous variable and to range from $m=0$ to
$m=1$ in the two phases.

Cells enter the proliferating phase with age $a=0$ and they are
committed to undergo cell division a time $\tau$ later, so the age
variable ranges from $a=0$ to $a=\tau$ in the proliferating phase.
We suppose that proliferating cells can be lost  by apoptosis with
a rate $\gamma$.

At the cytokinesis age, a cell divides and gives two daughter
cells, which enter immediatly the resting phase, with age $a=0$. A
cell can stay its entire life in the resting phase, so the age
variable ranges from $a=0$ to $a=+\infty$. The resting phase is a
quiescent stage in the cellular development. In this phase, cells
can either return to the proliferating phase at a rate $\beta $
and complete the cycle or die at a rate $\delta $ before ending
the cycle. According to a work of Sachs \cite{sachs}, we suppose
that the maturation of a cell and the density of resting cells at
a given maturity level determine the capacity of this cell for
entering the proliferating phase.

We denote by $p(t,m,a)$ and $n(t,m,a)$ respectively the population
densities in the proliferating and the resting phases at time $t$,
with age $a$ and maturity $m$. The conservation equations are
\begin{equation}
\frac{\partial p}{\partial t}+\frac{\partial p}{\partial
a}+\frac{\partial (V(m)p)}{\partial m} =-\gamma (m)p,
\label{prolif}
\end{equation}
\begin{equation}
\frac{\partial n}{\partial t}+\frac{\partial n}{\partial
a}+\frac{\partial (V(m)n)}{\partial m} =-\Big( \delta(m)+\beta
\big(m,N(t,m)\big)\Big) n,  \label{resting}
\end{equation}
where $V(m)$ is the maturation velocity and $N(t,m)$ is the
density of resting cells at time $t$ with a maturity level $m$,
defined by
\begin{displaymath}
N(t,m)=\displaystyle\int_{0}^{+\infty }n(t,m,a)da.
\end{displaymath}
We suppose that the function $V$ is continuously differentiable on
$\left[ 0,1\right] $, positive on $\left( 0,1\right] $ and
satisfies $V(0)=0$ and
\begin{equation}
\displaystyle\int_{0}^{m}\dfrac{ds}{V(s)}=+\infty , \qquad \text{
for }m\in \left( 0,1\right].  \label{(H.1).1}
\end{equation}
Since $\int_{m_{1}}^{m_{2}}\frac{ds}{V(s)}$, with $m_1<m_2$, is the time required for a cell with maturity $m_{1}$ to reach the maturity $m_{2}$, then Condition (\ref{(H.1).1}) means that a cell with very small maturity needs a long time to become mature.\\
For example, if
\begin{displaymath}
V(m)\underset{m\to 0}{\sim} \alpha m^p,\quad \textrm{ with }
\alpha>0 \textrm{ and } p\geq 1,
\end{displaymath}
then Condition (\ref{(H.1).1}) is satisfied.\\
We suppose, throughout this paper, that $\gamma$ and $\delta$ are
continuous and non-negative on $[0,1]$. The function $\beta$ is
supposed to be positive and continuous.

Equations (\ref{prolif}) and (\ref{resting}) are completed by
boundary conditions which represent the cellular flux between the
two phases. The first condition,
\begin{equation}
p(t,m,0)=\int_{0}^{+\infty }\beta
(m,N(t,m))n(t,m,a)da=\beta(m,N(t,m))N(t,m),  \label{bctauvariant1}
\end{equation}
describes the efflux of cells leaving the resting phase to the
proliferating one. Cells entering the proliferating phase with age
$0$ depend only on the population of the resting phase with a
given maturity level.

The second boundary condition determines the transfer of cells from the point of cytokinesis to the resting compartment.\\
We assume that a cell entering the proliferating phase with a
maturity $m\in[0,1]$ divides at age $\tau(m)>0$, and we require
that $\tau$ is a continuously differentiable and positive function
on $[0,1]$ such that
\begin{equation}
\tau^{\prime }(m)+\frac{1}{V(m)}>0, \qquad \textrm{ for }
m\in(0,1]. \label{taucondition}
\end{equation}
Since $V(0)=0$, this condition is always satisfied in a
neighborhood of the origin. If we suppose that the less mature
cells divide faster than more mature cells, that is, if we assume,
for example, that $\tau$ is an increasing function, then Condition
(\ref{taucondition}) is also satisfied.

\noindent If one consider a cell in the proliferating phase at
time $t$, with maturity $m\in(0,1]$, age $a$ and initial maturity
(that means at age $a=0$) $m_{0}$, then, naturally, we have
\begin{displaymath}
m_{0}\leq m \qquad \textrm{ and } \qquad
a=\int_{m_{0}}^{m}\frac{ds}{V(s)}\leq \tau(m_{0}).
\end{displaymath}
If $m$ is the maturity of the cell at the cytokinesis point, then
there exists a unique $\Theta(m)\in(0,m)$ (the maturity at the
point of commitment) such that
\begin{equation}\label{thetam}
\int_{\Theta(m)}^{m}\frac{ds}{V(s)}=\tau(\Theta(m)),
\end{equation}
because Condition (\ref{taucondition}) implies that the function
\begin{displaymath}
\widetilde{m} \rightarrow
\int_{\widetilde{m}}^{m}\frac{ds}{V(s)}-\tau(\widetilde{m})
\end{displaymath}
is continuous and strictly decreasing from $(0,m]$ into
$[-\tau(m),+\infty)$. Then, we can define a function $\Theta :
(0,1] \to (0,1]$, where $\Theta(m)$ satisfies (\ref{thetam}).

\noindent From a biological point of view, $\Theta(m)$ represents
the initial maturity of proliferating cells that divide at
maturity $m$ (at the point of cytokinesis). Then, from the
definition, the age of a cell with maturity $m$ at the point of
cytokinesis is $\tau(\Theta(m))$.

\noindent Remark that $\Theta$ is continuously differentiable on
$(0,1]$ and satisfies
\begin{displaymath}
0<\Theta (m)<m, \qquad \textrm{ for } m\in(0,1].
\end{displaymath}
This implies, in particular, that
\begin{equation} \label{propertytheta}
\lim_{m\to 0}\Theta(m)=0 \qquad \textrm{ and } \qquad \lim_{m\to
0}\int_{\Theta(m)}^{m}\frac{ds}{V(s)} =\tau(0)<+\infty.
\end{equation}

\noindent The property (\ref{propertytheta}) means that cells with
null maturity at the point of commitment keep a null maturity in
the proliferating phase.

\noindent The total number of proliferating cells at time $t$,
with maturity $m$, is given by
\begin{displaymath}
P(t,m)=\int_0^{\tau(\Theta(m))} p(t,m,a) da.
\end{displaymath}

\noindent We consider the characteristic curves $\chi :
(-\infty,0]\times[0,1] \to [0,1]$, solutions of the ordinary
differential equation
\begin{displaymath}
\left\{ \begin{array}{rcll}
\displaystyle\frac{d\chi}{ds}(s,m) & = & V(\chi(s,m)),& s \leq 0 \textrm{ and } m\in[0,1], \\
\chi(0,m)        & = & m. &
\end{array} \right.
\end{displaymath}
They represent the evolution of cells maturity to reach a maturity
$m$ at time $0$ from a time $s \leq 0$. They satisfy $\chi(s,0)=0$
and $\chi(s,m) \in (0,1]$ for $s\leq 0$ and $m \in (0,1]$.

\noindent It is not difficult to verify that, if $m\in[0,1]$, then
$\Theta (m)$ is the unique solution of the equation
\begin{equation} \label{propertyoftheta}
x=\chi(-\tau(x),m).
\end{equation}

\noindent At the end of the proliferating phase, a cell with a
maturity $m$ divides into two daughter cells with maturity $g(m)$.
We assume that $g:[0,1] \rightarrow [0,1]$ is a continuous and
strictly increasing function, continuously differentiable on
$[0,1)$ and such that $g(m)\leq m$ for $m\in [0,1]$. We also
assume, for technical reason and without loss of generality, that
\begin{displaymath}
\lim_{m\to 1}g^{\prime }(m)=+\infty.
\end{displaymath}
Then we can set
\begin{displaymath}
g^{-1}(m)=1, \qquad \textrm{ for } m>g(1).
\end{displaymath}
This means that the function $g^{-1}:[0,1] \rightarrow [0,1]$ is
continuously differentiable and satisfies
\begin{displaymath}
(g^{-1})^{\prime }(m)=0, \qquad \textrm{ for } m>g(1).
\end{displaymath}
Note that the maturity $m$ of the daughter cells just after
division is smaller than $g(1)$. Then, we must have
\begin{equation}
n(t,m,0)=0,\qquad \textrm{ for } m>g(1).  \label{For2}
\end{equation}
If a daughter cell has a maturity $m$ at birth, then the maturity
of its mother at the point of cytokinesis was $g^{-1}(m)$ and, at
the point of commitment, it was $\Theta (g^{-1}(m))$. We set
\begin{equation}
\Delta(m)=\Theta(g^{-1}(m)), \qquad \textrm{ for } m\in[0,1].
\label{equationdelta}
\end{equation}
From a biological point of view, $\Delta$ gives the link between
the maturity of a new born cell and the maturity of its mother at
the point of commitment. $\Delta : [0,1] \rightarrow [0,1] $ is
continuous, continuously differentiable on $(0,1]$, with
$\Delta(0)=0$. Moreover, $\Delta$ is strictly increasing on
$(0,g(1))$ with $\Theta(m)\leq \Delta(m)$ and
$\Delta(m)=\Theta(1)$ for $m\in[g(1),1]$.

\noindent Then, we can give the second boundary condition,
\begin{equation}
n(t,m,0)=2(g^{-1})^{\prime }(m)p\big(t,g^{-1}(m),\tau(\Delta
(m))\big), \quad \textrm{ for } t\geq 0 \textrm{ and } m\in[0,1].
\label{bctauvariant2}
\end{equation}
One can note that Expression (\ref{bctauvariant2}) includes also
Condition (\ref{For2}).

To complete the description of the model, we specify initial
conditions,
\begin{equation} \label{icp}
p(0,m,a)=\Gamma(m,a), \quad \textrm{ for }
(m,a)\in[0,1]\times[0,\tau_{max}],
\end{equation}
and
\begin{equation} \label{icn}
n(0,m,a)=\mu(m,a), \quad \textrm{ for }
(m,a)\in[0,1]\times[0,+\infty),
\end{equation}
where $\tau_{max}:=\max_{m\in[0,1]} \tau(m) >0$. $\Gamma $ and
$\mu $ are assumed to be continuous, and the function
\begin{equation} \label{gammamubarre}
\overline{\mu} : m\mapsto \int_{0}^{+\infty }\mu (m,a)da
\end{equation}
is supposed to be continuous on $[0,1]$.

\noindent We put
\begin{displaymath}
\xi(t,m):=\exp \left\{ -\int_{0}^{t}\Big(
\gamma\big(\chi(-s,m)\big)+V^{\prime }\big(\chi(-s,m)\big)\Big)
ds\right\},  \label{defksi}
\end{displaymath}
for $t\geq0$ and $m\in[0,1]$, and we define the sets
\begin{displaymath}
\Omega_{\Delta} := \Big\{ (m,t)\in[0,1]\times[0,+\infty) \ ; \
0\leq t \leq \tau(\Delta (m)) \Big\},
\end{displaymath}
and
\begin{displaymath}
\Omega_{\Theta} := \Big\{ (m,t)\in[0,1]\times[0,+\infty) \ ; \
0\leq t \leq \tau(\Theta (m)) \Big\}.
\end{displaymath}

\begin{proposition} \label{propositionequationsmodel}
Assume that the initial conditions $\mu$ and $\Gamma$ satisfy, for
$m\in[0,1]$,
\begin{equation} \label{condcont}
\Gamma(m,0)=\beta\big( m, \overline{\mu}(m)\big)
\overline{\mu}(m).
\end{equation}
Then, the total populations of proliferating and resting cells,
$P(t,m)$ and $N(t,m)$, satisfy, for $m\in[0,1]$ and $t\geq 0$,
\begin{equation}\label{equationP}
\begin{array}{l}
\displaystyle\frac{\partial }{\partial t}P(t,m)+\displaystyle\frac{\partial }{\partial m}(V(m)P(t,m))=-\gamma(m)P(t,m)+\beta \big(m, N(t,m)\big)N(t,m)\\
\quad \\
-\left\{ \setlength\arraycolsep{2pt} \begin{array}{ll}
\displaystyle \pi(m)\xi(t,m)\Gamma\Big( \chi\big(-t,m\big),\tau(\Theta(m))-t\Big), & \textrm{if } (m,t)\in\Omega_{\Theta},\\
\quad \\
\displaystyle \pi(m)\xi\big(\tau(\Theta(m)),m\big)\beta\Big(\Theta(m),N\big(t-\tau(\Theta(m)),\Theta(m)\big)\Big)\times\\
\quad \\
\quad\qquad\qquad\qquad\qquad\qquad\qquad
N\big(t-\tau(\Theta(m)),\Theta(m)\big), & \textrm{if }
(m,t)\notin\Omega_{\Theta},
\end{array} \right.
\end{array}
\end{equation}
\begin{equation}\label{equationN}
\begin{array}{l}
\displaystyle\frac{\partial }{\partial t}N(t,m)+\displaystyle\frac{\partial }{\partial m}(V(m)N(t,m))=-\Big(\delta (m)+\beta \big(m, N(t,m)\big)\Big)N(t,m)\\
\quad \\
+ \left\{ \begin{array}{ll}
\displaystyle 2(g^{-1})^{\prime }(m)\xi(t,g^{-1}(m))\Gamma\Big( \chi\big(-t,g^{-1}(m)\big),\tau(\Delta(m))-t\Big), & \textrm{if } (m,t)\in\Omega_{\Delta},\\
\quad \\
\displaystyle
\zeta(m)\beta\Big(\Delta(m),N\big(t-\tau(\Delta(m)),\Delta(m)\big)\Big)N\big(t-\tau(\Delta(m)),\Delta(m)\big),
& \textrm{if } (m,t)\notin\Omega_{\Delta},
\end{array} \right.
\end{array}
\end{equation}
and
\begin{eqnarray}
P(0,m)&=&\overline{\Gamma}(m):= \int_0^{\tau(\Theta (m))} \Gamma(m,a)da,\label{initialconditionP} \\
N(0,m)&=&\overline{\mu}(m),\label{initialconditionN}
\end{eqnarray}
with
\begin{displaymath}
\pi(m)=\frac{1}{1+V(\Theta(m))\tau^{\prime}(\Theta(m))},
\end{displaymath}
and
\begin{equation}\label{defzeta}
\zeta(m)=2(g^{-1})^{\prime
}(m)\xi\big(\tau(\Delta(m)),g^{-1}(m)\big).
\end{equation}
\end{proposition}

\vspace{1ex}

\begin{proof}
Using (\ref{icp}), (\ref{icn}) and the definitions of $P$ and $N$,
we obtain immediatly the equations (\ref{initialconditionP}) and
(\ref{initialconditionN}).

\noindent System (\ref{prolif})-(\ref{resting}) can be solved by
using the method of characteristics. First, we obtain the
following representation of solutions of Equation (\ref{prolif}),
\begin{displaymath}
p(t,m,a)=\left\{ \begin{array}{ll}
\xi(t,m)p(0,\chi(-t,m),a-t), &\textrm{for } 0\leq t<a, \\
\xi(a,m)p(t-a,\chi(-a,m),0), &\textrm{for } a\leq t.
\end{array} \right.  
\end{displaymath}
The initial condition (\ref{icp}) and the boundary condition
(\ref{bctauvariant1}) give
\begin{equation}
\begin{array}{l}
p(t,m,a)= \\
\quad \left\{ \begin{array}{ll}
\xi(t,m)\Gamma\big(\chi(-t,m),a-t\big), &\textrm{for } 0\leq t<a, \\
\quad & \quad \\
\xi(a,m)\beta\Big(\chi(-a,m),N\big(t-a,\chi(-a,m)\big)\Big)N\big(t-a,\chi(-a,m)\big),
&\textrm{for } a\leq t.
\end{array} \right.  \label{Caract1}
\end{array}
\end{equation}

\noindent Let $m\in[0,1]$ be given. By integrating Equation
(\ref{prolif}) with respect to the age, between $0$ and
$\tau(\Theta(m))$, we obtain
\begin{displaymath}
\frac{\partial }{\partial
t}P(t,m)+\int_{0}^{\tau(\Theta(m))}\frac{\partial }{\partial
m}(V(m)p(t,m,a)) da =
-\gamma(m)P(t,m)+p(t,m,0)-p\big(t,m,\tau(\Theta(m))\big).
\end{displaymath}
One can note that
\begin{displaymath}
\tau^{\prime}(\Theta(m))\Theta^{\prime}(m)V(m)-1=-\pi(m).
\end{displaymath}
Since
\begin{displaymath}
\displaystyle\frac{\partial }{\partial
m}(V(m)P(t,m))=\displaystyle\int_{0}^{\tau(\Theta(m))}\frac{\partial
}{\partial m}\big(V(m)p(t,m,a)) da
+\tau^{\prime}(\Theta(m))\Theta^{\prime}(m)V(m)p\big(t,m,\tau(\Theta(m))\big),
\end{displaymath}
and
\begin{displaymath}
\begin{array}{l}
p(t,m,\tau(\Theta(m)))= \\
\\
\quad \left\{ \begin{array}{ll}
\xi(t,m)\Gamma(\chi(-t,m),\tau(\Theta(m))-t), &\textrm{if } 0\leq t<\tau(\Theta(m)), \\
\quad & \quad \\
\xi(\tau(\Theta(m)),m)\beta\Big(\chi(-\tau(\Theta(m)),m),N\big(t-\tau(\Theta(m)),\chi(-\tau(\Theta(m)),m)\big)\Big)\times \\
\quad & \quad \\
\quad\quad\quad\quad\quad\quad\quad\quad\quad\quad\quad\quad\quad\quad\quad\quad
N\big(t-\tau(\Theta(m)),\chi(-\tau(\Theta(m)),m)\big), &\textrm{if
} \tau(\Theta(m))\leq t,
\end{array} \right.
\end{array}
\end{displaymath}
then, using (\ref{bctauvariant1}) and (\ref{propertyoftheta}), we
obtain Equation (\ref{equationP}).

\noindent Thanks to (\ref{condcont}) and by using the continuity
of $\overline{\mu}$, we show that $\lim_{a\to +\infty}
n(t,m,a)=0$.

\noindent So, by integrating Equation (\ref{resting}) with respect
to the age, between $0$ and $+\infty$, it follows that
\begin{displaymath}
\frac{\partial }{\partial t}N(t,m)+\frac{\partial }{\partial
m}(V(m)N(t,m))
=-\Big(\delta(m)+\beta(m,N(t,m))\Big)N(t,m)+n(t,m,0).
\end{displaymath}
From the equations (\ref{propertyoftheta}) and
(\ref{equationdelta}), we deduce that
\begin{displaymath}
\Delta (m)=\chi\big(-\tau(\Delta(m)),g^{-1}(m)\big).
\end{displaymath}
Hence, from Equations (\ref{bctauvariant2}) and (\ref{Caract1}),
we obtain
\begin{displaymath}
n(t,m,0)=\left\{ \begin{array}{ll}
\displaystyle 2(g^{-1})^{\prime }(m)\xi(t,g^{-1}(m))\Gamma\Big( \chi\big(-t,g^{-1}(m)\big),\tau(\Delta(m))-t\Big), & \textrm{if } (m,t)\in\Omega_{\Delta},\\
\quad \\
\displaystyle
\zeta(m)\beta\Big(\Delta(m),N\big(t-\tau(\Delta(m)),\Delta(m)\big)\Big)N\big(t-\tau(\Delta(m)),\Delta(m)\big),
& \textrm{if } (m,t)\notin\Omega_{\Delta}.
\end{array} \right.
\end{displaymath}
Equation (\ref{equationN}) follows immediatly.

\noindent Finally, we can remark that, if $N$ is continuous,
Condition (\ref{condcont}) implies that the mappings $(m,t)\mapsto
F\Big(t,m,N\big(t-\tau(\Delta(m)),\Delta(m)\big)\Big)$ and
$(m,t)\mapsto
G\Big(t,m,N\big(t-\tau(\Theta(m)),\Theta(m)\big)\Big)$, with
$F:[0,+\infty)\times[0,1]\times\mathbb{R}\to \mathbb{R}$ and
$G:[0,+\infty)\times[0,1]\times\mathbb{R}\to \mathbb{R}$ given by
\begin{equation} \label{F}
F(t,m,x)=\left\{ \begin{array}{ll}
\displaystyle 2(g^{-1})^{\prime }(m)\xi(t,g^{-1}(m))\Gamma\Big( \chi\big(-t,g^{-1}(m)\big),\tau(\Delta(m))-t\Big), & \textrm{if } (m,t)\in\Omega_{\Delta},\\
\quad \\
\displaystyle \zeta(m)\beta\big(\Delta(m),x\big)x, & \textrm{if }
(m,t)\notin\Omega_{\Delta},
\end{array} \right.
\end{equation}
and
\begin{equation} \label{G}
G(t,m,x)= \left\{ \setlength\arraycolsep{2pt} \begin{array}{ll}
\displaystyle \pi(m)\xi(t,m)\Gamma\Big( \chi\big(-t,m\big),\tau(\Theta(m))-t\Big), & \textrm{if } (m,t)\in\Omega_{\Theta},\\
\quad \\
\displaystyle
\pi(m)\xi\big(\tau(\Theta(m)),m\big)\beta\big(\Theta(m),x\big)x, &
\textrm{if } (m,t)\notin\Omega_{\Theta},
\end{array} \right.
\end{equation}
are continuous.

\noindent This completes the proof.
\end{proof}

\noindent One can remark that the solutions of Equations
(\ref{equationN}) and (\ref{initialconditionN}) do not depend on
the proliferating population. We extend $N$ by setting
\begin{equation} \label{prolongementN}
N(t,m)=\overline{\mu}(m), \quad \textrm{ for } t\in[-\tau_{max},
0] \textrm{ and } m\in[0,1].
\end{equation}
This extension does not influence our system. However, it will be
useful in the following.

\section{Local existence and global continuation} \label{localexistence}

In this section, we are interested in proving the local existence
of an integrated solution of Problem
(\ref{equationP})-(\ref{initialconditionN}). First, we consider an
integrated formulation of Problem
(\ref{equationP})-(\ref{initialconditionN}). We denote by $C[0,1]$
the space of continuous functions on $[0,1]$, endowed with the
supremum norm $||.||$, defined by
\begin{equation}\label{supnorm}
||v||=\sup_{m\in[0,1]} |v(m)|, \qquad \textrm{ for } v\in C[0,1].
\end{equation}

\noindent Let us consider the unbounded closed linear operator $A
: D(A)\subset C[0,1]\to C[0,1]$ defined by
\begin{displaymath}
D(A)=\big\{ u\in C[ 0,1] \ ; u \textrm{ differentiable on }(0,1],
u^{\prime }\in C(0,1], \ \lim_{x\rightarrow 0}V(x)u^{\prime }(x)=0
\big\}
\end{displaymath}
and
\begin{displaymath}
Au(x)=\left\{ \begin{array}{ll}
-(\delta(x)+V^{\prime }(x))u(x)-V(x)u^{\prime }(x), & \quad \textrm{if } x\in(0,1], \\
-(\delta(0)+V^{\prime }(0))u(0), & \quad \textrm{if }x=0.
\end{array}\right.
\end{displaymath}

\noindent Then, we have the following proposition, which
characterise the operator $(A,D(A))$.

\begin{proposition} \label{propsemigroupe}
The operator $A$ is the infinitesimal generator of the strongly
continuous semigroup $(T(t))_{t\geq 0}$ defined on $C[0,1]$ by
\begin{displaymath}
(T(t)\psi)(x)=K(t,x)\psi(\chi(-t,x)), \quad \textrm{ for } \psi\in
C[0,1], t\geq 0 \textrm{ and } x\in[0,1],
\end{displaymath}
where
\begin{displaymath}
K(t,x)=\exp \left\{ -\int_{0}^{t}\Big(
\delta\big(\chi(-s,x)\big)+V^{\prime }\big(\chi(-s,x)\big)\Big)ds
\right\}.
\end{displaymath}
\end{proposition}

\begin{proof}
The proof is similar to the proof of Proposition 2.4 in
\cite{webb1996}.
\end{proof}

\noindent We denote by $C(\Omega_{\Theta})$ the space of
continuous function on $\Omega_{\Theta}$, endowed with the norm
$\| . \|_{\Omega_{\Theta}}$, defined by
\begin{displaymath}
\| \Upsilon \|_{\Omega_{\Theta}}:=
\sup_{(m,a)\in\Omega_{\Theta}}|\Upsilon(m,a)|, \quad \textrm{ for
} \Upsilon \in C(\Omega_{\Theta}).
\end{displaymath}

\noindent Now, we can consider an integrated formulation of
Problem (\ref{equationP})-(\ref{initialconditionN}), given by the
variation of constant formula associated to the $C_0$-semigroup
$(T(t))_{t\geq 0}$. That is the following definition.

\begin{definition}
Let $\Gamma\in C(\Omega_{\Theta})$ and $\mu$ be a function such
that $\overline{\mu}\in C[0,1]$, with $\overline{\mu}$ given by
(\ref{gammamubarre}). An \emph{integrated solution} of Problem
(\ref{equationP})-(\ref{initialconditionN}) is a continuous
solution of the system
\begin{equation} \label{integratedformulationN}
\begin{array}{l}
N(t,m)=K(t,m)\overline{\mu}\big(\chi(-t,m)\big)\\
 \\
-\displaystyle\int_0^t \!\! K(t-s,m)\beta\Big(\chi(-(t-s),m),N\big(s,\chi(-(t-s),m)\big)\Big)N\big(s,\chi(-(t-s),m)\big) ds\\
 \\
+\displaystyle\int_0^t \!\! K(t-s,m)
F\Big(s,\chi(-(t-s),m),N\big(s-\tau(\Delta(\chi(-(t-s),m))),\Delta(\chi(-(t-s),m))\big)\Big)
ds,
\end{array}
\end{equation}
and
\begin{equation} \label{integratedformulationP}
\begin{array}{l}
P(t,m)=\xi(t,m)\overline{\Gamma}\big(\chi(-t,m)\big)\\
 \\
+\displaystyle\int_0^t \!\! \xi(t-s,m)\beta\Big(\chi(-(t-s),m),N\big(s,\chi(-(t-s),m)\big)\Big)N\big(s,\chi(-(t-s),m)\big) ds\\
 \\
-\displaystyle\int_0^t \!\! \xi(t-s,m)
G\Big(s,\chi(-(t-s),m),N\big(s-\tau(\Theta(\chi(-(t-s),m))),\Theta(\chi(-(t-s),m))\big)\Big)
ds,
\end{array}
\end{equation}
for $t\geq 0$ and $m\in[0,1]$, where $F$ and $G$ are given by
(\ref{F}) and (\ref{G}) and $\overline{\Gamma}$ is given by
(\ref{initialconditionP}).
\end{definition}

\noindent The extension given by (\ref{prolongementN}) allows the
second integrals, in the expressions
(\ref{integratedformulationN}) and (\ref{integratedformulationP}),
to be well defined.

\noindent In order to obtain a result of local existence for the
solutions of System
(\ref{integratedformulationN})-(\ref{integratedformulationP}), we
first focus on  Equation (\ref{integratedformulationN}). We show,
in the next theorem, that Equation (\ref{integratedformulationN})
has a unique local solution, which depends continuously on the
initial conditions.

\begin{theorem} \label{theoremlocalexistence}
Assume that the mapping $x\mapsto x\beta(m,x)$ is locally
Lipschitz continuous for all $m\in[0,1]$, that is, for all $r>0$,
there exists $L(r)\geq 0$ such that
\begin{displaymath} 
|x\beta(m,x)-y\beta(m,y)|\leq L(r)|x-y|, \quad \textrm{ if }
|x|<r, |y|<r \textrm{ and } m\in[0,1].
\end{displaymath}
If $\Gamma\in C(\Omega_{\Theta})$ and $\mu$ is a function such
that $\overline{\mu}\in C[0,1]$, then, there exists $T_{max}>0$
such that Equation (\ref{integratedformulationN}) has a unique
continuous solution $N^{\overline{\mu},\Gamma}$ defined on a
maximal domain $[0,T_{max})\times[0,1]$, and either
\begin{displaymath}
T_{max}=+\infty \qquad \textrm{ or } \qquad  \limsup_{t\to
T_{max}^-} \|N^{\overline{\mu},\Gamma}(t,.)\|=+\infty.
\end{displaymath}
Furthermore, $N^{\overline{\mu},\Gamma}(t,.)$ is a continuous
function of $\overline{\mu}$ and $\Gamma$, in the sense that, if
$t\in(0,T_{max})$, $\overline{\mu}_1\in C[0,1]$ and $\Gamma_1\in
C(\Omega_{\Theta})$, then there exist a continuous positive
function $C:[0,+\infty)\to\mathbb{R}$ and a constant
$\varepsilon>0$ such that, for $\overline{\mu}_2 \in C[0,1]$ and
$\Gamma_2 \in C(\Omega_{\Theta})$ such that
$N^{\overline{\mu}_2,\Gamma_2}$ is defined on $[0,t]\times[0,1]$
and
\begin{displaymath}
\|\overline{\mu}_1 - \overline{\mu}_2 \|<\varepsilon \quad
\textrm{ and } \quad  \|\Gamma_1 - \Gamma_2
\|_{\Omega_{\Theta}}<\varepsilon,
\end{displaymath}
we get
\begin{displaymath}
\|
N^{\overline{\mu}_1,\Gamma_1}(s,.)-N^{\overline{\mu}_2,\Gamma_2}(s,.)
\| \leq C(t)\Big( \|\overline{\mu}_1 - \overline{\mu}_2 \|+
\|\Gamma_1 - \Gamma_2 \|_{\Omega_{\Theta}}\Big), \quad \textrm{
for } s\in[0,t].
\end{displaymath}
\end{theorem}

\begin{proof}
We put
\begin{displaymath}
r=\| \overline{\mu} \|+1.
\end{displaymath}
Let $T>0$ be fixed. We consider the following set,
\begin{displaymath}
X(\overline{\mu})=\bigg\{ N\in C([0,T]\times[0,1])\ ; \
N(0,.)=\overline{\mu} \textrm{ on } [0,1] \ \textrm{ and }
\sup_{(t,m)\in[0,T]\times[0,1]} |N(t,m)-\overline{\mu}(m) | \leq 1
\bigg\},
\end{displaymath}
where $C([0,T]\times[0,1])$ is endowed with the uniform norm. $X(\overline{\mu})$ is a non-empty  closed convex subset of $C([0,T]\times[0,1])$.\\
We define the operator $H : C([0,T]\times[0,1])\to
C([0,T]\times[0,1])$ by
\begin{displaymath}
\begin{array}{l}
H(N)(t,m)=K(t,m)\overline{\mu}\big(\chi(-t,m)\big)\\
 \\
-\displaystyle\int_0^t K(t-s,m)\beta\Big(\chi(-(t-s),m),N\big(s,\chi(-(t-s),m)\big)\Big)N\big(s,\chi(-(t-s),m)\big) ds\\
 \\
+\displaystyle\int_0^t K(t-s,m)
F\Big(s,\chi(-(t-s),m),N\big(s-\tau(\Delta(\chi(-(t-s),m))),\Delta(\chi(-(t-s),m))\big)\Big)
ds.
\end{array}
\end{displaymath}
$H$ is continuous in $C([0,T]\times[0,1])$.
Our objective is to show that $H$ is a contraction from $X(\overline{\mu})$ into itself.\\
Let $N\in X(\overline{\mu})$. It is clear that
$H(N)(0,.)=\overline{\mu}$. On the other hand, we have, for
$(t,m)\in[0,T]\times[0,1]$,
\begin{displaymath}
\begin{array}{l}
\big|H(N)(t,m)-\displaystyle\overline{\mu}(m)\big|\leq\big|K(t,m)\overline{\mu}\big(\chi(-t,m)\big)-\overline{\mu}(m)\big|\\
 \\
+\displaystyle\bigg|\int_0^t K(t-s,m)\beta\Big(\chi(-(t-s),m),N\big(s,\chi(-(t-s),m)\big)\Big)N\big(s,\chi(-(t-s),m)\big) ds\bigg|\\
 \\
+\displaystyle\bigg|\int_0^t K(t-s,m)
F\Big(s,\chi(-(t-s),m),N\big(s-\tau(\Delta(\chi(-(t-s),m))),\Delta(\chi(-(t-s),m))\big)\Big)
ds\bigg|.
\end{array}
\end{displaymath}
Since $K$ is continuous on $[0,T]\times[0,1]$, then there exists
$\widetilde{K}\geq 0$ such that
\begin{displaymath}
|K(t,m)|\leq \widetilde{K}, \qquad \textrm{ for }
(t,m)\in[0,T]\times[0,1].
\end{displaymath}
Since $N\in X(\overline{\mu})$, then
\begin{displaymath}
|N(t,m)|\leq 1+\| \overline{\mu} \|=r.
\end{displaymath}
This implies that there exists $M:=\max\{\widetilde{M}, \| \zeta
\| rL(r) \}\geq 0$ such that, for $(t,m)\in[0,T]\times[0,1]$ and
$s\in[0,t]$,
\begin{displaymath}
\Big|
F\Big(s,\chi(-(t-s),m),N\big(s-\tau(\Delta(\chi(-(t-s),m))),\Delta(\chi(-(t-s),m))\big)\Big)
\Big| \leq M,
\end{displaymath}
where
\begin{displaymath}
\widetilde{M}:=\sup_{(m,t)\in\Omega_{\Delta}}
\Big|2(g^{-1})^{\prime }(m)\xi(t,g^{-1}(m))\Gamma\Big(
\chi\big(-t,g^{-1}(m)\big),\tau(\Delta(m))-t\Big)\Big|,
\end{displaymath}
and $\zeta$ is given by (\ref{defzeta}). Hence, we obtain that
\begin{displaymath}
\big|H(N)(t,m)-\overline{\mu}(m)\big|\leq
\big|K(t,m)\overline{\mu}\big(\chi(-t,m)\big)-\overline{\mu}(m)\big|+\widetilde{K}(rL(r)+M)t.
\end{displaymath}
Let us recall that $K(0,m)=1$, $\chi(0,m)=m$ and $(t,m)\mapsto
K(t,m)\overline{\mu}\big(\chi(-t,m)\big)$ is continuous. Then, we
can choose $T>0$ such that
\begin{equation} \label{inequality}
\displaystyle\sup_{(t,m)\in[0,T]\times[0,1]} \bigg\{
\big|K(t,m)\overline{\mu}\big(\chi(-t,m)\big)-\overline{\mu}(m)\big|+\widetilde{K}(rL(r)+M)t
\bigg\} < 1.
\end{equation}
Consequently,
\begin{displaymath}
\big|H(N)(t,m)-\overline{\mu}(m)\big|\leq 1,\quad \textrm{ for }
(t,m)\in[0,T]\times[0,1],
\end{displaymath}
and $H(X(\overline{\mu}))\subset X(\overline{\mu})$.\\
Now, we show that $H$ is a contraction on $X(\overline{\mu})$.\\
Let $N_{1}\in X(\overline{\mu})$ and $N_{2}\in X(\overline{\mu})$.
Then,
\begin{displaymath}
\setlength\arraycolsep{2pt}
\begin{array}{l}
|H(N_{1})(t,m)-H(N_{2})(t,m)|\\
\quad \\
\leq \bigg|\displaystyle\int_0^t K(t-s,m)\bigg[\beta\Big(\chi(-(t-s),m),N_1\big(s,\chi(-(t-s),m)\big)\Big)N_1\big(s,\chi(-(t-s),m)\big)\\
\quad \\
\quad \qquad \qquad \qquad \qquad -\beta\Big(\chi(-(t-s),m),N_2\big(s,\chi(-(t-s),m)\big)\Big)N_2\big(s,\chi(-(t-s),m)\big)\bigg] ds \bigg|\\
\quad \\
+ \bigg|\displaystyle\int_0^t K(t-s,m)\bigg[ F\Big(s,\chi(-(t-s),m),N_1\big(s-\tau(\Delta(\chi(-(t-s),m))),\Delta(\chi(-(t-s),m))\big)\Big)\\
\quad \\
\quad \quad \qquad \qquad -F\Big(s,\chi(-(t-s),m),N_2\big(s-\tau(\Delta(\chi(-(t-s),m))),\Delta(\chi(-(t-s),m))\big)\Big) \bigg] ds\bigg|,\\
\quad \\
\leq \widetilde{K}(1+\| \zeta \|)L(r)T
\displaystyle\sup_{(t,m)\in[0,T]\times[0,1]}|N_1(t,m)-N_2(t,m)|.
\end{array}
\end{displaymath}
Since $r\geq 1$ and $\| \zeta \|rL(r)\leq M$, then Condition
(\ref{inequality}) implies that
\begin{displaymath}
\widetilde{K}(1+\| \zeta \|) L(r)T \leq \widetilde{K}(1+\| \zeta
\|)rL(r)T \leq \widetilde{K}(rL(r)+M)T < 1.
\end{displaymath}
Hence, $H$ is a contraction from $X(\overline{\mu})$ into itself.
\noindent Therefore, there exists a unique $N\in X(\overline{\mu})$ such that $N$ satisfies Equation (\ref{integratedformulationN}) for $(t,m)\in[0,T]\times[0,1]$.\\
Let $N$ be the maximal solution of Equation (\ref{integratedformulationN}), defined on its maximal domain $[0,T_{max})\times[0,1]$.\\
Assume that
\begin{displaymath}
T_{max}<+\infty \qquad \textrm{ and } \qquad \limsup_{t\to
T^{-}_{max}} \|N(t,.)\|<+\infty.
\end{displaymath}
Then, there exists $r>0$ such that
\begin{displaymath}
\|N(t,.)\|<r, \quad \textrm{ for } t\in[0,T_{max}).
\end{displaymath}
Let $(t,m)\in[0,T_{max})\times[0,1]$ and $h>0$ such that
$t+h\in[0,T_{max})$. Then,
\begin{displaymath}
\begin{array}{l}
|N(t+h,m)-N(t,m)|\\
\quad\\
\leq \Big| K(t+h,m)\overline{\mu}\big(\chi(-(t+h),m)\big)-K(t,m)\overline{\mu}\big(\chi(-t,m)\big) \Big|\\
\quad\\
+\bigg|\displaystyle\int_0^{t+h} K(t+h-s,m)\beta\Big(\chi(-(t+h-s),m),N\big(s,\chi(-(t+h-s),m)\big)\Big)N\big(s,\chi(-(t+h-s),m)\big) ds \\
\quad\\
\quad\quad\quad\quad \ \  -\displaystyle\int_0^t K(t-s,m)\beta\Big(\chi(-(t-s),m),N\big(s,\chi(-(t-s),m)\big)\Big)N\big(s,\chi(-(t-s),m)\big) ds \bigg|\\
\quad\\
+\bigg|\displaystyle\int_0^{t+h} K(t+h-s,m)\times\\
\quad\\
\quad\quad F\Big(s,\chi(-(t+h-s),m),N\big(s-\tau(\Delta(\chi(-(t+h-s),m))),\Delta(\chi(-(t+h-s),m))\big)\Big) ds\\
\quad\\
\quad -\displaystyle\int_0^t K(t-s,m)F\Big(s,\chi(-(t-s),m),N\big(s-\tau(\Delta(\chi(-(t-s),m))),\Delta(\chi(-(t-s),m))\big)\Big) ds\bigg|,\\
\quad\\
\leq \Big|
K(t+h,m)\overline{\mu}\big(\chi(-(t+h),m)\big)-K(t,m)\overline{\mu}\big(\chi(-t,m)\big)
\Big| + |I_1| + |I_2|.
\end{array}
\end{displaymath}
Since
\begin{displaymath}
\begin{array}{l}
I_1:=\displaystyle\int_0^{t+h} K(t+h-s,m)\beta\Big(\chi(-(t+h-s),m),N\big(s,\chi(-(t+h-s),m)\big)\Big)N\big(s,\chi(-(t+h-s),m)\big) ds\\
\quad \\
 -\displaystyle\int_0^t K(t-s,m)\beta\Big(\chi(-(t-s),m),N\big(s,\chi(-(t-s),m)\big)\Big)N\big(s,\chi(-(t-s),m)\big) ds, \\
\quad \\
=\displaystyle\int_{-h}^{t} K(t-s,m)\beta\Big(\chi(-(t-s),m),N\big(s+h,\chi(-(t-s),m)\big)\Big)N\big(s+h,\chi(-(t-s),m)\big) ds\\
\quad \\
 -\displaystyle\int_0^t K(t-s,m)\beta\Big(\chi(-(t-s),m),N\big(s,\chi(-(t-s),m)\big)\Big)N\big(s,\chi(-(t-s),m)\big) ds, \\
\quad \\
=\displaystyle\int_{-h}^{0} K(t-s,m)\beta\Big(\chi(-(t-s),m),N\big(s+h,\chi(-(t-s),m)\big)\Big)N\big(s+h,\chi(-(t-s),m)\big) ds\\
\quad \\
 +\displaystyle\int_0^t K(t-s,m)\bigg[ \beta\Big(\chi(-(t-s),m),N\big(s+h,\chi(-(t-s),m)\big)\Big)N\big(s+h,\chi(-(t-s),m)\big)\\
\quad \\
\qquad\quad\quad\quad\quad\quad\quad\quad-
\beta\Big(\chi(-(t-s),m),N\big(s,\chi(-(t-s),m)\big)\Big)N\big(s,\chi(-(t-s),m)\big)\bigg]ds,
\end{array}
\end{displaymath}
then
\begin{displaymath}
|I_1| \leq
\widetilde{K}rL(r)h+\widetilde{K}L(r)\displaystyle\int_0^t
\Big|N\big(s+h,\chi(-(t-s),m)\big)-N\big(s,\chi(-(t-s),m)\big)
\Big|ds.
\end{displaymath}
By using the same reasonning, we also obtain that
\begin{displaymath}
\begin{array}{rcl}
|I_2| &\leq& \widetilde{K}Mh+\varrho_1(h) +\widetilde{K}\| \zeta \|L(r)\displaystyle\int_0^t \Big|N\big(s+h-\tau(\Delta(\chi(-(t-s),m))),\Delta(\chi(-(t-s),m))\big)\\
      &  & \quad\qquad\qquad\qquad\qquad\qquad\qquad -N\big(s-\tau(\Delta(\chi(-(t-s),m))),\Delta(\chi(-(t-s),m))\big) \Big|ds,
\end{array}
\end{displaymath}
with $\varrho_1(h)$ independent of $(t,m)$ and such that
$\lim_{h\to 0} \varrho_1(h)=0$.

\noindent One can remark that, if $t_{1}>0$, $t_{2}>0$ and
$m\in[0,1]$, then
\begin{displaymath}
K(t_{1}+t_{2},m)=K(t_{1},m)K(t_{2},\chi(-t_{1},m))=K(t_{2},m)K(t_{1},\chi(-t_{2},m)).
\end{displaymath}
This yields to
\begin{displaymath}
\begin{array}{l}
|N(t+h,m)-N(t,m)|\\
\quad \\
 \leq \widetilde{K}\Big(\varrho_2(h)+rL(r)h+Mh \Big)+ \varrho_1(h)+\widetilde{K}L(r)\displaystyle\int_0^t \| N(s+h, .)-N(s,.) \| ds\\
\quad\\
+\widetilde{K}\| \zeta \|L(r)\displaystyle\int_0^t \Big\|
N\Big(s+h-\tau\big(\Delta(\chi(-(t-s),m))\big),
.\Big)-N\Big(s-\tau\big(\Delta(\chi(-(t-s),m))\big),.\Big)
\Big\|ds,
\end{array}
\end{displaymath}
where
\begin{displaymath}
\varrho_2(h)=\sup_{(t,m)\in[0,T_{max})\times[0,1]} \Big|
K\big(h,\chi(-t,m)\big)\overline{\mu}\big(\chi(-(t+h),m)\big)-\overline{\mu}\big(\chi(-t,m)\big)
\Big| \underset{h\to 0}{\to} 0.
\end{displaymath}
We set
\begin{displaymath}
\varrho(h)=\widetilde{K}\Big(\varrho_2(h)+rL(r)h+Mh
\Big)+\varrho_1(h).
\end{displaymath}
Hence, for $(t,m)\in[0,T_{max})\times[0,1]$, $h>0$ and
$\theta\in[-\tau_{max},0]$ such that $t+h\in[0,T_{max})$ and
$t+\theta\geq 0$,
\begin{displaymath}
\begin{array}{l}
|N(t+h+\theta,m)-N(t+\theta,m)| \leq \varrho(h)+\widetilde{K}L(r)\displaystyle\int_0^{t+\theta} \| N(s+h, .)-N(s,.) \| ds\\
+\widetilde{K}\| \zeta \|L(r)\displaystyle\int_0^{t+\theta} \|
N\big(s+h-\tau(\Delta(\chi(-(t-s),m))),
.\big)-N\big(s-\tau(\Delta(\chi(-(t-s),m))),.\big) \|ds.
\end{array}
\end{displaymath}
On the other hand,
\begin{displaymath}
N(t,m)=\overline{\mu}(m), \quad \textrm{ for } t\in[-\tau_{max},0]
\textrm{ and } m\in[0,1].
\end{displaymath}
This implies that
\begin{displaymath}
\displaystyle\sup_{\theta\in[-\tau_{max},0]}\|
N(t+h+\theta,.)-N(t+\theta,.) \| \leq \varrho(h)+\widetilde{K}(\|
\zeta \|+1)L(r)\displaystyle\int_0^{t}
\sup_{\theta\in[-\tau_{max},0]}\| N(s+h+\theta, .)-N(s+\theta,.)
\| ds.
\end{displaymath}
By using the Gronwall's inequality, it follows that
\begin{displaymath}
\displaystyle\sup_{\theta\in[-\tau_{max},0]}\|
N(t+h+\theta,.)-N(t+\theta,.) \| \leq
\varrho(h)e^{\widetilde{K}(\| \zeta \|+1)L(r)T_{max}}.
\end{displaymath}
Hence,
\begin{displaymath}
\| N(t+h,.)-N(t,.) \| \leq \varrho(h) e^{\widetilde{K}(\| \zeta
\|+1)L(r)T_{max}},
\end{displaymath}
with
\begin{displaymath}
\lim_{h\to 0} \varrho(h) = 0.
\end{displaymath}
Using the same reasonning, we can show a same result for $h<0$.\\
It follows immediatly that
\begin{displaymath}
\lim_{t\to T^{-}_{max}} N(t,.) \textrm{ exists}.
\end{displaymath}
This implies that $N$ can be extended continuously to $T_{max}$,
which contradicts the maximality of $[0,T_{max})$. Hence,
\begin{displaymath}
T_{max}=+\infty \qquad \textrm{ or } \qquad \limsup_{t\to
T^{-}_{max}} \|N(t,.)\|=+\infty.
\end{displaymath}
Now, we denote by $N^{\overline{\mu}_1,\Gamma_1}$ the solution of
Equation (\ref{integratedformulationN}) for the initial data
$\overline{\mu}_1$ and $\Gamma_1$, defined on its maximal domain
$[0,T_{max})\times[0,1]$. Let $t\in(0,T_{max})$ be fixed. We put
\begin{displaymath}
r(t)=\sup_{s\in[0,t]} \| N^{\overline{\mu}_1,\Gamma_1}(s,.)\|
\qquad \textrm{ and } \qquad R(t)=1+r(t).
\end{displaymath}
We set
\begin{displaymath}
\widetilde{C}:=\sup_{(m,t)\in\Omega_{\Theta}} |2(g^{-1})^{\prime
}(m)\xi(t,g^{-1}(m))|,
\end{displaymath}
and
\begin{displaymath}
C(t)=\widetilde{K}\max\{1,\widetilde{C}\}e^{L(R(t))\widetilde{K}(1+\|
\zeta \|)t}.
\end{displaymath}
Let $0<\varepsilon<1$ be such that $2C(t)\varepsilon\in(0,1)$ and
let $\overline{\mu}_2 \in C[0,1]$ and $\Gamma_2 \in
C(\Omega_{\Theta})$ such that $N^{\overline{\mu}_2,\Gamma_2}$ is
defined on $[0,t]\times[0,1]$ and
\begin{displaymath}
\|\overline{\mu}_1 - \overline{\mu}_2 \|<\varepsilon \quad
\textrm{ and } \quad \|\Gamma_1 - \Gamma_2
\|_{\Omega_{\Theta}}<\varepsilon.
\end{displaymath}
Then, there exists $s>0$ such that
$\|N^{\overline{\mu}_2,\Gamma_2}(\sigma,.)\|\leq R(t)$ for all
$\sigma\in[0,s]$. Let
\begin{displaymath}
t_0 = \sup \bigg\{ s>0 \ ; \
\|N^{\overline{\mu}_2,\Gamma_2}(\sigma,.)\|\leq R(t), \ \textrm{
for } \sigma\in[0,s] \bigg\}.
\end{displaymath}
If we suppose that $t_0 <t$, then $N^{\overline{\mu}_2,\Gamma_2}$
is defined on $[0,t_0]\times[0,1]$ and satisfies, for
$s\in[0,t_0]$ and $m\in[0,1]$,
\begin{displaymath}
\begin{array}{rcl}
|N^{\overline{\mu}_1,\Gamma_1}(s,m)-N^{\overline{\mu}_2,\Gamma_2}(s,m)|&\leq& \widetilde{K}\Big(\|\overline{\mu}_1 - \overline{\mu}_2 \|+\widetilde{C}\|\Gamma_1 - \Gamma_2 \|_{\Omega_{\Theta}}\Big)\\
    &+&\widetilde{K}(1+\| \zeta \|)L(R(t))\displaystyle\int_{0}^{s} \|N^{\overline{\mu}_1,\Gamma_1}(\sigma,.)-N^{\overline{\mu}_2,\Gamma_2}(\sigma,.)\| d\sigma.
\end{array}
\end{displaymath}
Therefore, using the Gronwall's Inequality, we obtain
\begin{equation} \label{condcontinuability}
|N^{\overline{\mu}_1,\Gamma_1}(s,m)-N^{\overline{\mu}_2,\Gamma_2}(s,m)|\leq
C(t)\Big(\|\overline{\mu}_1-\overline{\mu}_2 \|+\|\Gamma_1 -
\Gamma_2 \|_{\Omega_{\Theta}}\Big), \quad \textrm{ for }
s\in[0,t_0].
\end{equation}
This implies, in particular, that
\begin{displaymath}
\|N^{\overline{\mu}_2,\Gamma_2}(s,.)\|\leq 2C(t)\varepsilon +
\|N^{\overline{\mu}_1,\Gamma_1}(s,.)\| < 1+r(t) = R(t), \quad
\textrm{ for } s\in[0,t_0].
\end{displaymath}
This contradicts the definition of $t_0$. Hence, $t_0\geq t$.\\
That means that (\ref{condcontinuability}) is satisfied for each
$t\in[0,T_{max})$. Then, we deduce the continuous dependence of
the solution with the initial data and the proof is complete.
\end{proof}

One can remark that Condition (\ref{condcont}) is not needed to
prove Theorem \ref{theoremlocalexistence}.

By using Theorem \ref{theoremlocalexistence}, we can deduce the
following result, which deals with the existence of solutions of
Problem
(\ref{integratedformulationN})-(\ref{integratedformulationP}).
\begin{corollary} \label{localexistenceP}
Under the assumptions of Theorem \ref{theoremlocalexistence}, System (\ref{integratedformulationN})-(\ref{integratedformulationP}) has a unique continuous maximal solution $(N,P)$, defined on $[0,T_{max})\times[0,1]$.
\end{corollary}

\begin{proof}
Under the assumptions of Theorem \ref{theoremlocalexistence}, Equation (\ref{integratedformulationN}) has a unique continuous maximal solution $N$ defined on  $[0,T_{max})\times[0,1]$. Then, we easily obtain the existence and uniqueness of a solution of Equation (\ref{integratedformulationP}) on  $[0,T_{max})\times[0,1]$.
\end{proof}

We can use the results of local existence, given by Theorem
\ref{theoremlocalexistence} and Corollary \ref{localexistenceP},
to investigate the global existence of the solutions of System
(\ref{integratedformulationN})-(\ref{integratedformulationP}).
This is done in the next theorem.

\begin{theorem}\label{theoremglobalexistence}
Under the assumptions of Theorem \ref{theoremlocalexistence} and
the assumption that the mapping $x\mapsto \beta(m,x)$ is uniformly
bounded, the unique solution of System
(\ref{integratedformulationN})-(\ref{integratedformulationP}) is
global, that means, it is defined for all $t\geq 0$.
\end{theorem}

\begin{proof}
We assume that there exists $\widetilde{\beta}\geq 0$ such that
\begin{displaymath}
|\beta(m,x)|\leq \widetilde{\beta}, \qquad \textrm{ for all }
m\in[0,1] \textrm{ and } x\in\mathbb{R}.
\end{displaymath}
Then, for $t\in[0,T_{max})$ and $m\in[0,1]$, we get
\begin{displaymath}
\begin{array}{rcl}
|N(t,m)|&\leq&\widetilde{K} \| \overline{\mu} \|+\widetilde{K}\widetilde{\beta}\displaystyle\int_0^t \| N(s,.) \| ds\\
    &    & \quad \\
    & +  &\widetilde{K}\widetilde{M}t+\widetilde{K}\widetilde{\beta}\| \zeta \|\displaystyle\int_0^t \| N\big(s-\tau(\Delta(\chi(-(t-s),m))),. \big) \| ds.
\end{array}
\end{displaymath}
Let $\theta\in[-\tau_{max},0]$ be such that $t+\theta\geq 0$.
Then,
\begin{displaymath}
|N(t+\theta,m)|\leq \widetilde{K}(\| \overline{\mu}
\|+\widetilde{M}t)+\widetilde{K}\widetilde{\beta}(\| \zeta
\|+1)\displaystyle\int_0^t
\sup_{\overline{\theta}\in[-\tau_{max},0]}\|
N(s+\overline{\theta},.) \| ds.
\end{displaymath}
On the other hand, if $-\tau_{max} \leq t+\theta\leq 0$, then
\begin{displaymath}
|N(t+\theta,m)|\leq \widetilde{K}\| \overline{\mu} \|,
\end{displaymath}
because $\widetilde{K}\geq 1$. Hence,
\begin{displaymath}
\sup_{\theta\in[-\tau_{max},0]} \| N(t+\theta,.) \|\leq
\widetilde{K}(\| \overline{\mu}
\|+\widetilde{M}t)+\widetilde{K}\widetilde{\beta}(\| \zeta
\|+1)\displaystyle\int_0^t \sup_{\theta\in[-\tau_{max},0]}\|
N(s+\theta,.) \| ds.
\end{displaymath}
By using the Gronwall's inequality, we deduce that
\begin{displaymath}
\| N(t,.) \| \leq  \sup_{\theta\in[-\tau_{max},0]} \|
N(t+\theta,.) \| \leq \widetilde{K}(\| \overline{\mu}
\|+\widetilde{M}t)e^{\widetilde{K}\widetilde{\beta}(\| \zeta
\|+1)t}:= f(t).
\end{displaymath}
Since $f$ is continuous on $[0,T_{max}]$, then
\begin{displaymath}
\limsup_{t\to T_{max}^-} \| N(t,.) \| <+\infty.
\end{displaymath}
We deduce that $T_{max}=+\infty$ and the solution of
(\ref{integratedformulationN}) is global.

\noindent Finally, we easily obtain that the unique maximal
solution of Equation (\ref{integratedformulationP}) is also
global.
\end{proof}

\begin{corollary}\label{globalsolution}
Assume that the mapping $x\mapsto x\beta(m,x)$ is Lipschitz
continuous for all $m\in[0,1]$. Then, for $\Gamma\in
C(\Omega_{\Theta})$ and $\mu$ such that $\overline{\mu}\in
C[0,1]$,  the unique solution of System
(\ref{integratedformulationN})-(\ref{integratedformulationP}) is
global.
\end{corollary}

We have studied, throughout Theorems \ref{theoremlocalexistence}
and \ref{theoremglobalexistence}, the local and global existence
of the solutions of (\ref{integratedformulationN}) and of Problem
(\ref{integratedformulationN})-(\ref{integratedformulationP}).
Before we investigate the positivity of these solutions in the
next section, we can ask for regularity results for these
solutions. This is presented in the following remark.

\begin{remark}
Under classical assumptions on the function $\beta$ and the
initial data $\Gamma$ and $\mu$, we can obtain regularity results
for the solution of System
(\ref{integratedformulationN})-(\ref{integratedformulationP}).
This may be done by using the same idea as in the work of Travis
and Webb \cite{traviswebb}.
\end{remark}

\section{Positivity of solutions} \label{sectionpositivityregularity}

Since we study a biological population, it is necessary to obtain
the positivity of the solutions of System
(\ref{integratedformulationN})-(\ref{integratedformulationP}) to
ensure that the model is well-posed. First, we focus our study on
the solutions of Equation (\ref{integratedformulationN}). We use a
method given by Webb \cite{webb1985} in 1985 and developped by
Kato \cite{kato2000} to obtain the positivity of the solutions of
Equation (\ref{integratedformulationN}).

Let $\Gamma\in C(\Omega_{\Theta})$ and $\mu$ be a function such
that $\overline{\mu}\in C[0,1]$, $\overline{\mu}$ given by
(\ref{gammamubarre}). Let $T>0$ be fixed. We consider the family
of operators $H^{a} : C([0,T)\times[0,1]) \to C([0,T)\times[0,1])$
defined, for $a\in C[0,1]$, $(t,m)\in[0,T)\times[0,1]$ and $N\in
C([0,T)\times[0,1])$, by
\begin{equation} \label{Ha}
\begin{array}{l}
H^{a}(N)(t,m) = e^{-ta(\chi(-t,m))}K(t,m)\overline{\mu}\big(\chi(-t,m)\big)\\
 \\
+\displaystyle\int_0^t e^{-(t-s)a(\chi(-t,m))}K(t-s,m)\bigg[a\Big(\chi\big(-(t-s),m\big)\Big)N\Big(s,\chi\big(-(t-s),m\big)\Big)\\
 \\
-\beta\bigg(\chi\big(-(t-s),m\big),N\Big(s,\chi\big(-(t-s),m\big)\Big)\bigg)N\Big(s,\chi\big(-(t-s),m\big)\Big)\\
 \\
+F\Big(s,\chi\big(-(t-s),m\big),N\big(s-\tau(\Delta(\chi(-(t-s),m))),\Delta(\chi(-(t-s),m))\big)\Big)
\bigg] ds.
\end{array}
\end{equation}
Let $a\in C[0,1]$, $(t,m)\in(0,T)\times[0,1]$ and $N\in
C([0,T)\times[0,1])$ be fixed. We consider the mapping $w^{a} :
[0,t] \to \mathbb{R}$ defined, for $s\in[0,t]$, by
\begin{displaymath}
w^{a}(s)=H^{a}(N)\Big(s,\chi\big(-(t-s),m\big)\Big).
\end{displaymath}
Then, we can prove the following lemmas.

\begin{lemma}\label{lemma1}
The function $w^{a}$ is differentiable on $[0,t]$ and satisfies
\begin{equation} \label{wderivee}
\left\{ \begin{array}{rcl}
\displaystyle \frac{d}{ds}w^{a}(s)&=&-\nu^a(s)w^{a}(s)+f^a(s), \quad \textrm{ for } s\in[0,t],\\
w^a(0)&=&\overline{\mu}\big(\chi(-t,m)\big),
\end{array} \right.
\end{equation}
with
\begin{displaymath}
\setlength\arraycolsep{2pt}
\begin{array}{rcl}
f^a(s)&=&a\Big(\chi\big(-(t-s), m\big)\Big)N\Big(s,\chi\big(-(t-s),m\big) \Big)\\
 &&\\
&-& N\Big(s, \chi\big(-(t-s),m\big) \Big)\beta\bigg( \chi\big(-(t-s),m\big) ,N\Big( s, \chi\big(-(t-s),m\big) \Big)\bigg) \\
 &&\\
&+& F\Bigg(s,\chi\big(-(t-s),m\big) ,
N\bigg(s-\tau\Big(\Delta\Big( \chi\big(-(t-s),m\big)
\Big)\Big),\Delta\Big( \chi\big(-(t-s),m\big)\Big)\bigg)\Bigg),
\end{array}
\end{displaymath}
and
\begin{displaymath}
\nu^a(s)=a\big(\chi(-t,m)\big)+\delta\Big(\chi\big(-(t-s),m\big)\Big)+V^{\prime}\Big(\chi\big(-(t-s),m\big)\Big).
\end{displaymath}
\end{lemma}

\begin{proof}
First, one can remark that, if $t_{1}$ and $t_{2}$ are positive,
then
\begin{displaymath}
\chi\big(-(t_{1}+t_{2}),m \big)=\chi\big(-t_{1}, \chi\big(-t_{2},m
\big) \big).
\end{displaymath}
It follows that, for $s\in[0,t]$ and $h\neq 0$ such that
$s+h\in[0,t]$, we have
\begin{displaymath}
\frac{1}{h}\big(w^{a}(s+h)-w^{a}(s)\big)=W_{1,h}+W_{2,h},
\end{displaymath}
with
\begin{displaymath}
\setlength\arraycolsep{2pt}
\begin{array}{l}
W_{1,h}=\displaystyle\frac{1}{h}\bigg[
e^{-ha(\chi(-t,m))}K\Big(h,\chi\big(-(t-s-h),m\big)\Big)-1
\Bigg]e^{-sa(\chi(-t,m))}K\Big(s,\chi\big(-(t-s),m\big)\Big)\overline{\mu}\big(\chi(-t,m)\big),
\end{array}
\end{displaymath}
and
\begin{displaymath}
\setlength\arraycolsep{2pt}
\begin{array}{rcl}
W_{2,h}&=&\displaystyle\frac{1}{h}\displaystyle\int_0^{s+h}\!\!\! e^{-(s+h-\sigma)a(\chi(-t, m))}K\Big(s+h-\sigma,\chi\big(-(t-s-h),m\big)\Big)\times\\
&& \\
&&\Bigg[a\Big(\chi\big(-(t-\sigma), m\big)\Big)N\Big(\sigma,\chi\big(-(t-\sigma), m\big)\Big)\\
&& \\
&&- N\Big(\sigma,\chi\big(-(t-\sigma), m\big)\Big)\beta\bigg(\chi\big(-(t-\sigma), m \big),N\Big(\sigma,\chi\big(-(t-\sigma), m\big)\Big) \bigg) \\
 \\
&&+F\Bigg(\sigma,\chi\big(-(t-\sigma), m \big), N\bigg(\sigma-\tau\Big(\Delta\Big(\chi\big(-(t-\sigma), m\big)\Big)\Big),\Delta\Big(\chi\big(-(t-\sigma), m \big)\Big)\bigg)\Bigg)\Bigg] d\sigma\\
&& \\
&-&\displaystyle\frac{1}{h}\displaystyle\int_0^{s}\!\!\! e^{-(s-\sigma)a(\chi(-t, m))}K\Big(s-\sigma,\chi\big(-(t-s),m\big)\Big)\times\\
&& \\
&&\Bigg[a\Big(\chi\big(-(t-\sigma), m\big)\Big) N\Big(\sigma,\chi\big(-(t-\sigma), m\big)\Big)\\
&& \\
&&- N\Big(\sigma,\chi\big(-(t-\sigma), m\big)\Big)\beta\bigg(\chi\big(-(t-\sigma), m \big),N\Big(\sigma,\chi\big(-(t-\sigma), m\big)\Big) \bigg) \\
&& \\
&&+F\Bigg(\sigma,\chi\big(-(t-\sigma), m \big),
N\bigg(\sigma-\tau\Big(\Delta\Big(\chi\big(-(t-\sigma),
m\big)\Big)\Big),\Delta\Big(\chi\big(-(t-\sigma), m
\big)\Big)\bigg)\Bigg) \Bigg] d\sigma.
\end{array}
\end{displaymath}
Hence, we easily obtain that
\begin{displaymath}
\lim_{h\to 0} W_{1,h} = -\nu^a(s)e^{-sa(\chi(-t,m))}
K\Big(s,\chi\big(-(t-s),m\big)\Big)\overline{\mu}\big(\chi(-t,m)\big).
\end{displaymath}
By the same way, we have
\begin{displaymath}
\setlength\arraycolsep{2pt}
\begin{array}{l}
W_{2,h}=\displaystyle\frac{1}{h}\bigg[ e^{-ha(\chi(-t,m))}K\Big(h,\chi\big(-(t-s-h),m\big)\Big)-1 \bigg]\times\\
 \\
\displaystyle\int_0^{s} e^{-(s-\sigma)a(\chi(-t,m))}K\Big(s-\sigma,\chi\big(-(t-s),m\big)\Big)\Bigg[a\Big(\chi\big(-(t-\sigma), m\big)\Big)N\Big(\sigma,\chi\big(-(t-\sigma),m\big) \Big)\\
 \\
- N\Big( \sigma, \chi\big(-(t-\sigma),m\big) \Big)\beta\bigg( \chi\big(-(t-\sigma),m\big) ,N\Big( \sigma, \chi\big(-(t-\sigma),m\big) \Big)\bigg) \\
 \\
+ F\Bigg(\sigma,\chi\big(-(t-\sigma),m\big) , N\bigg(\sigma-\tau\Big(\Delta\Big( \chi\big(-(t-\sigma),m\big) \Big)\Big),\Delta\Big( \chi\big(-(t-\sigma),m\big)\Big)\bigg)\Bigg) \Bigg] d\sigma \\
 \\
+\displaystyle\frac{1}{h}e^{-ha(\chi(-t,m))}K\Big(h,\chi\big(-(t-s-h),m\big)\Big)\displaystyle\int_s^{s+h} e^{-(s-\sigma)a(\chi(-t,m))}K\Big(s-\sigma,\chi\big(-(t-s),m\big)\Big)\times\\
 \\
\Bigg[a\Big(\chi\big(-(t-\sigma), m\big)\Big)N\Big(\sigma,\chi\big(-(t-\sigma),m\big) \Big)\\
 \\
- N\Big( \sigma, \chi\big(-(t-\sigma),m\big) \Big)\beta\bigg( \chi\big(-(t-\sigma),m\big) ,N\Big( \sigma, \chi\big(-(t-\sigma),m\big) \Big)\bigg) \\
 \\
+ F\Bigg(\sigma,\chi\big(-(t-\sigma),m\big) ,
N\bigg(\sigma-\tau\Big(\Delta\Big( \chi\big(-(t-\sigma),m\big)
\Big)\Big),\Delta\Big(
\chi\big(-(t-\sigma),m\big)\Big)\bigg)\Bigg) \Bigg] d\sigma.
\end{array}
\end{displaymath}
Therefore,
\begin{displaymath}
\begin{array}{rcl}
\displaystyle\lim_{h\to 0} W_{2,h}&=& -\nu^a(s)\displaystyle\int_0^{s} e^{-(s-\sigma)a(\chi(-t,m))}K\Big(s-\sigma,\chi\big(-(t-s),m\big)\Big)\times\\
&& \\
&&\Bigg[a\Big(\chi\big(-(t-\sigma), m\big)\Big)N\Big(\sigma,\chi\big(-(t-\sigma),m\big) \Big)\\
&& \\
&&- N\Big( \sigma, \chi\big(-(t-\sigma),m\big) \Big)\beta\bigg( \chi\big(-(t-\sigma),m\big) ,N\Big( \sigma, \chi\big(-(t-\sigma),m\big) \Big)\bigg) \\
&& \\
&&+ F\Bigg(\sigma,\chi\big(-(t-\sigma),m\big) , N\bigg(\sigma-\tau\Big(\Delta\Big( \chi\big(-(t-\sigma),m\big) \Big)\Big),\Delta\Big( \chi\big(-(t-\sigma),m\big)\Big)\bigg)\Bigg) \Bigg] d\sigma \\
&& \\
&&+a\Big(\chi\big(-(t-s), m\big)\Big)N\Big(s,\chi\big(-(t-s),m\big) \Big)\\
&& \\
&&- N\Big(s, \chi\big(-(t-s),m\big) \Big)\beta\bigg( \chi\big(-(t-s),m\big) ,N\Big( s, \chi\big(-(t-s),m\big) \Big)\bigg) \\
&& \\
&&+ F\Bigg(s,\chi\big(-(t-s),m\big) ,
N\bigg(s-\tau\Big(\Delta\Big( \chi\big(-(t-s),m\big)
\Big)\Big),\Delta\Big( \chi\big(-(t-s),m\big)\Big)\bigg)\Bigg).
\end{array}
\end{displaymath}
Hence, $\lim_{h\to 0}\frac{1}{h}(w^{a}(s+h)-w^{a}(s))$ exists and is equal to $-\nu^a(s)w^{a}(s)+f^a(s)$.\\
We obtain Equation (\ref{wderivee}) and the proof of the lemma is
complete.
\end{proof}

\begin{lemma} \label{lemma2}
Let $a,b\in C[0,1]$. Then,
\begin{displaymath}
\begin{array}{rcl}
H^{b}(N)(t,m)&=&H^{a}(N)(t,m)+\displaystyle\int_{0}^{t}  e^{-(t-s)b(\chi(-t,m))}K(t-s,m)\times\\
    & & \Big(b\big(\chi(-(t-s),m)\big)-a\big(\chi(-(t-s),m)\big)\Big)\bigg[ N\Big(s,\chi\big(-(t-s),m\big)\Big)-w^{a}(s) \bigg] ds.
\end{array}
\end{displaymath}
\end{lemma}

\begin{proof}
From Equation (\ref{wderivee}), it follows that
\begin{displaymath}
\setlength\arraycolsep{2pt}
\begin{array}{rcl}
\displaystyle\frac{d}{ds}\Big[w^{b}(s)-w^{a}(s)\Big]&=&-\nu^b(s)\big(w^{b}(s)-w^{a}(s)\big)\\
& & +
\Big(b\big(\chi(-(t-s),m)\big)-a\big(\chi(-(t-s),m)\big)\Big)\bigg[
N\Big(s,\chi\big(-(t-s),m\big)\Big)-w^{a}(s)  \bigg].
\end{array}
\end{displaymath}
Since $w^{b}(0)=w^{a}(0)=\overline{\mu}\big(\chi(-t,m)\big)$, then
\begin{displaymath}
\begin{array}{l}
\big(w^{b}(t)-w^{a}(t)\big)e^{\int_{0}^{t}\nu^b(\sigma)d\sigma}=\\
\displaystyle\int_{0}^{t} e^{\int_{0}^{s}\nu^b(\sigma)d\sigma}
\Big(b\big(\chi(-(t-s),m)\big)-a\big(\chi(-(t-s),m)\big)\Big)
\bigg[ N\Big(s,\chi\big(-(t-s),m\big)\Big)-w^{a}(s)  \bigg]ds.
\end{array}
\end{displaymath}
Therefore,
\begin{displaymath}
\begin{array}{rcl}
w^{b}(t)-w^{a}(t)&=&\displaystyle\int_{0}^{t}  e^{-(t-s)b(\chi(-t,m))}K(t-s,m)\times\\
&
&\Big(b\big(\chi(-(t-s),m)\big)-a\big(\chi(-(t-s),m)\big)\Big)\bigg[
N\Big(s,\chi\big(-(t-s),m\big)\Big)-w^{a}(s)  \bigg] ds.
\end{array}
\end{displaymath}
This ends the proof.
\end{proof}

\begin{lemma}\label{lemma3}
Assume that there exist $a\in C[0,1]$ and $N\in
C([0,T)\times[0,1])$ such that $H^a(N)=N$. Then, for all $b\in
C[0,1]$, $H^b(N)=N$.
\end{lemma}

\begin{proof}
This is a simple consequence of Lemma \ref{lemma2}.
\end{proof}

Now, we show, under a cellular regulation hypothesis, that the
solutions of
(\ref{integratedformulationN})-(\ref{integratedformulationP}) are
positive. We use the three previous lemmas to show the following
theorem.

\begin{theorem}  \label{theorempositivity}
Assume that the mapping $x\mapsto \beta(m,x)$ satisfies the
cellular regulation hypothesis,
\begin{equation} \label{betapropertyfeedback}
\big(\beta(m,x)-\beta(m,0)\big)x\leq 0, \quad  \textrm{ for }
m\in[0,1] \textrm{ and } x\in\mathbb{R}.
\end{equation}
Let $\Gamma\in C(\Omega_{\Theta})$ and $\mu$ be a continuous
function such that $\overline{\mu}\in C[0,1]$. If
$\overline{\mu}\geq 0$ and $\Gamma\geq 0$, then the solutions of
Equation (\ref{integratedformulationN}) are non-negative on their
domains.
\end{theorem}

\begin{proof}
Let $N$ be a solution of Equation (\ref{integratedformulationN}) defined on $[0,T)\times[0,1]$, where $T>0$.\\
We consider the operator $H^{a}$ defined, for $a\in C[0,1]$, by
(\ref{Ha}). It is obvious that
\begin{displaymath}
H^{0}(N)=N,
\end{displaymath}
where $0$ denotes the null mapping of $C[0,1]$.

\noindent Then, from Lemma \ref{lemma3}, for all $a\in C[0,1]$, $H^{a}(N)=N$. Hence, if there exists $a\in C[0,1]$ such that $H^{a}(N)(t,m)\geq 0$ for $(t,m)\in[0,T)\times[0,1]$, then we shall obtain that $N(t,m)\geq 0$, for $(t,m)\in[0,T)\times[0,1]$.\\
\noindent We set
\begin{displaymath}
a(m):= \beta(m,0), \qquad \textrm{ for all } m\in[0,1].
\end{displaymath}
Then, $a\in C[0,1]$ and Condition (\ref{betapropertyfeedback})
implies that
\begin{displaymath}
a(m)N(t,m)-N(t,m)\beta(m,N(t,m))\geq 0, \quad \textrm{ for }
(t,m)\in[0,T)\times[0,1].
\end{displaymath}
Hence,
\begin{displaymath}
\begin{array}{l}
\displaystyle\int_0^t e^{-(t-s)a(\chi(-t,m))}K(t-s,m)\bigg[a\Big(\chi\big(-(t-s),m\big)\Big)N\Big(s,\chi\big(-(t-s),m\big)\Big)\\
-\beta\bigg(\chi\big(-(t-s),m\big),N\Big(s,\chi\big(-(t-s),m\big)\Big)\bigg)N\Big(s,\chi\big(-(t-s),m\big)\Big)\bigg]
ds
\end{array}
\end{displaymath}
is non-negative. We have to show by steps that, for
$(t,m)\in[0,T)\times[0,1]$,
\begin{displaymath}
\int_0^te^{-(t-s)a(\chi(-t,m))} K(t-s,m)
F\Big(s,\chi\big(-(t-s),m\big),N\big(s-\tau(\Delta(\chi(-(t-s),m))),\Delta(\chi(-(t-s),m))\big)\Big)
ds
\end{displaymath}
is non-negative.\\
First, it is clear that $H^{a}(N)$ is non-negative on $\Omega_{\Delta}$.\\
We set
\begin{displaymath}
\tau_{\Delta} := \min_{m\in[0,1]} \tau(\Delta(m)).
\end{displaymath}
Since $\tau$ is positive on $[0,1]$, then $\tau_{\Delta}>0$.\\
Let $m\in[0,1]$. Suppose that
$t\in[\tau(\Delta(m)),\tau(\Delta(m))+\tau_{\Delta}]$. Then,
\begin{displaymath}
0\leq s-\tau\big(\Delta\big(\chi\big(-(t-s),m\big)\big)\big) \leq
t-\tau_{\Delta} \leq \tau(\Delta(m)), \quad \textrm{ for }
s\in[\tau(\Delta(m)),t].
\end{displaymath}
Since $F\geq 0$ for $(m,t)\in\Omega_{\Delta}$, this yields that
\begin{displaymath}
\displaystyle\int_0^t e^{-(t-s)a(\chi(-t,m))} K(t-s,m)
F\Big(s,\chi(-(t-s),m),N\big(s-\tau(\Delta(\chi(-(t-s),m))),\Delta(\chi(-(t-s),m))\big)\Big)
ds\geq 0.
\end{displaymath}
Hence, $H^{a}(N)(t,m)\geq 0$ for $m\in[0,1]$ and
$t\in[0,\tau(\Delta(m))+\tau_{\Delta}]$. This implies that
\begin{displaymath}
N(t,m)\geq 0, \qquad \textrm{ for } m\in[0,1] \textrm{ and }
t\in[0,\tau(\Delta(m))+\tau_{\Delta}].
\end{displaymath}
By steps, we show that
\begin{displaymath}
N(t,m)\geq 0, \qquad \textrm{ for } m\in[0,1] \textrm{ and }
t\in[0,\tau(\Delta(m))+n\tau_{\Delta}]\textrm{ with }
n\in\mathbb{N}.
\end{displaymath}
There exists $n\in\mathbb{N}$ such that
\begin{displaymath}
\tau(\Delta(m))+n\tau_{\Delta}<T
\leq\tau(\Delta(m))+(n+1)\tau_{\Delta}.
\end{displaymath}
Consequently, $N\geq 0$ on its domain $[0,T)\times[0,1]$ and the
proof is complete.
\end{proof}

We deduce immediatly the following corollary.

\begin{corollary}
Under the assumptions of Theorem \ref{theorempositivity}, the
solutions of Problem
(\ref{integratedformulationN})-(\ref{integratedformulationP}) are
non-negative.
\end{corollary}

\section{Local and global stability}\label{sectionlocalstability}

We study, in this section, the local and global stability of the
trivial solution of Problem
(\ref{integratedformulationN})-(\ref{integratedformulationP}). We
first focus on Equation (\ref{integratedformulationN}).

\noindent We set
\begin{displaymath}
\widetilde{\delta}:=\inf_{m\in[0,1]}
\Big(\delta(m)+V^{\prime}(m)\Big), \quad
\kappa:=\sup_{m\in[0,1]}|(g^{-1})^{\prime}(m)| \quad \textrm{and}
\quad \widetilde{\gamma}:=\inf_{m\in[0,1]}
\Big(\gamma(m)+V^{\prime}(m)\Big).
\end{displaymath}
We assume that $\widetilde{\gamma}\geq 0$ and that the function
$x\mapsto x\beta(m,x)$ satisfies a Lipshitz condition in a
neighborhood of zero, that is, there exist $\varepsilon>0$ and
$L\geq 0$ such that
\begin{equation}\label{Lipshitzconditionneighborhoodzero}
|x\beta(m,x)-y\beta(m,y)|\leq L|x-y|, \quad \textrm{ if }
|x|<\varepsilon, |y|<\varepsilon \textrm{ and } m\in[0,1].
\end{equation}

\noindent Throughout this section, $N^{\overline{\mu},\Gamma}$
denotes a global solution of Equation
(\ref{integratedformulationN}) associated with the initial data
$\mu$ and $\Gamma$.

\noindent In order to obtain the local stability of the trivial
solution, we first prove a result of invariance for the solutions.

\begin{proposition}\label{invariance}
Assume that Condition (\ref{Lipshitzconditionneighborhoodzero}) is
satisfied and that
\begin{equation}\label{condinvariance}
L(1+2\kappa)<\widetilde{\delta}.
\end{equation}
Let $\Gamma\in C(\Omega_{\Theta})$ and $\mu$ be a continuous
function such that $\overline{\mu}\in C[0,1]$. If
$\|\overline{\mu}\|\leq \varepsilon$ and $\| \Gamma
\|_{\Omega_{\Theta}}\leq \varepsilon L$, where $\varepsilon$ and
$L$ are given by (\ref{Lipshitzconditionneighborhoodzero}), then
\begin{displaymath}
|N^{\overline{\mu},\Gamma}(t,m)|\leq \varepsilon, \quad
\textrm{for } t\geq0 \textrm{ and } m\in[0,1].
\end{displaymath}
\end{proposition}

\begin{proof}
We first notice that, under the assumption (\ref{condinvariance}),
$\widetilde{\delta}>0$. We consider the sequence defined by
\begin{displaymath}
N_0(t,m)=K(t,m)\overline{\mu}\big(\chi(-t,m)\big),
\end{displaymath}
and
\begin{displaymath}
N_n(t,m)=N_0(t,m)-I(N_{n-1})(t,m)+J(N_{n-1})(t,m),
\end{displaymath}
for $t\geq 0$, $m\in[0,1]$ and $n\in\mathbb{N}^*$, with
\begin{displaymath}
I(N)(t,m)=\displaystyle\int_0^t \!\!
K(t-s,m)\beta\Big(\chi(-(t-s),m),N\big(s,\chi(-(t-s),m)\big)\Big)N\big(s,\chi(-(t-s),m)\big)
ds,
\end{displaymath}
and
\begin{displaymath}
J(N)(t,m)=\displaystyle\int_0^t \!\! K(t-s,m)
F\Big(s,\chi(-(t-s),m),N\big(s-\tau(\Delta(\chi(-(t-s),m))),\Delta(\chi(-(t-s),m))\big)\Big)
ds.
\end{displaymath}
Then,
\begin{displaymath}
|N_0(t,m)| \leq e^{-\widetilde{\delta}t}\|\overline{\mu}\|\leq
e^{-\widetilde{\delta}t}\varepsilon \leq\varepsilon.
\end{displaymath}
By induction, we show that
\begin{equation}\label{induction}
|N_n(t,m)| \leq \varepsilon, \qquad \textrm{ for } t\geq0, \
m\in[0,1] \textrm{ and } n\in\mathbb{N}^*.
\end{equation}
We assume that $|N_n(t,m)| \leq \varepsilon$.

\noindent Let $(t,m)\in[0,+\infty)\times[0,1]$ be fixed.

\noindent One has to notice that, if $0\leq s\leq
\tau(\Delta(\chi(-(t-s),m)))$, then
\begin{displaymath}
\begin{array}{l}
|F\Big(s,\chi(-(t-s),m),N_n\big(s-\tau(\Delta(\chi(-(t-s),m))),\Delta(\chi(-(t-s),m))\big)\Big)|\\
\quad\\
\leq 2\kappa e^{-\widetilde{\gamma}s}\| \Gamma \|_{\Omega_{\Theta}},\\
\quad\\
\leq 2\kappa\varepsilon L,
\end{array}
\end{displaymath}
and, if $s>\tau(\Delta(\chi(-(t-s),m)))$, then
\begin{displaymath}
\begin{array}{l}
|F\Big(s,\chi(-(t-s),m),N_n\big(s-\tau(\Delta(\chi(-(t-s),m))),\Delta(\chi(-(t-s),m))\big)\Big)|\\
\quad\\
\leq \| \zeta \|\bigg|\beta\Big(\Delta(\chi(-(t-s),m)),N_n\big(s-\tau(\Delta(\chi(-(t-s),m))),\Delta(\chi(-(t-s),m))\big)\Big)\times\\
\quad\\
\quad\quad N_n\big(s-\tau(\Delta(\chi(-(t-s),m))),\Delta(\chi(-(t-s),m))\big)\bigg|,\\
\quad\\
\leq \| \zeta \|\varepsilon L,
\end{array}
\end{displaymath}
where $\| \zeta \|$ is given by (\ref{supnorm}). Since $\| \zeta
\|\leq 2\kappa$, we obtain that
\begin{displaymath}
|J(N_{n})(t,m)|\leq 2\kappa\varepsilon L\displaystyle\int_0^t \!\!
e^{-\widetilde{\delta}(t-s)}ds,
\end{displaymath}
and
\begin{displaymath}
\begin{array}{rcl}
|N_{n+1}(t,m)|&\leq& |N_0(t,m)|+|I(N_{n})(t,m)|+|J(N_{n})(t,m)|,\\
&&\\
    &\leq& \varepsilon e^{-\widetilde{\delta}t} +\varepsilon L(1+2\kappa)\displaystyle\int_0^t e^{-\widetilde{\delta}(t-s)}ds,\\
&&\\
    &\leq& \varepsilon e^{-\widetilde{\delta}t} +\varepsilon \displaystyle\frac{L}{\widetilde{\delta}}(1+2\kappa)(1-e^{-\widetilde{\delta}t}).
\end{array}
\end{displaymath}
By assumption,
\begin{displaymath}
\frac{L}{\widetilde{\delta}}(1+2\kappa)<1.
\end{displaymath}
Consequently,
\begin{displaymath}
|N_{n+1}(t,m)|\leq
\varepsilon\Big(e^{-\widetilde{\delta}t}\Big(1-\frac{L}{\widetilde{\delta}}(1+2\kappa)\Big)+\frac{L}{\widetilde{\delta}}(1+2\kappa)\Big)
\leq \varepsilon.
\end{displaymath}
We conclude that (\ref{induction}) is true for $n\in\mathbb{N}^*$.

\noindent By remarking that the sequence
$(N_n)_{n\in\mathbb{N}^*}$ converges to
$N^{\overline{\mu},\Gamma}$, then we obtain that
\begin{displaymath}
|N^{\overline{\mu},\Gamma}(t,m)|\leq \varepsilon, \quad
\textrm{for } t\geq0 \textrm{ and } m\in[0,1],
\end{displaymath}
which ends the proof.
\end{proof}

Before we prove the local stability of Equation
(\ref{integratedformulationN}), we rewrite the integrated solution
of (\ref{integratedformulationN}) by using a variation of constant
formula, for $t\geq \tau_{max}$ and $m\in[0,1]$. We obtain that
\begin{displaymath} 
\begin{array}{l}
N(t,m)=K(t-\tau_{max},m)N\Big(\tau_{max},\chi\big(-(t-\tau_{max}),m\big)\Big)\\
 \\
-\displaystyle\int_{\tau_{max}}^t \!\! K(t-s,m)\beta\Big(\chi(-(t-s),m),N\big(s,\chi(-(t-s),m)\big)\Big)N\big(s,\chi(-(t-s),m)\big) ds\\
 \\
+\displaystyle\int_{\tau_{max}}^t \!\! K(t-s,m) \zeta(\chi(-(t-s),m))N\big(s-\tau(\Delta(\chi(-(t-s),m))),\Delta(\chi(-(t-s),m))\big)\times\\
 \\
\beta\Big(\Delta(\chi(-(t-s),m)),N\big(s-\tau(\Delta(\chi(-(t-s),m))),\Delta(\chi(-(t-s),m))\big)\Big)ds.
\end{array}
\end{displaymath}
Hence, we can show the following theorem, which deals with the
local stability of the trivial solution of Equation
(\ref{integratedformulationN}).

\begin{theorem}\label{theoremlocalstabilityN}
Under the assumptions of Proposition \ref{invariance}, the trivial
solution of Equation (\ref{integratedformulationN}) is locally
exponentially stable, that is, there exist $\overline{t}>0$,
$\epsilon>0$, $c\geq0$ and $d\geq0$ such that, if
$\|\overline{\mu}\|<\epsilon$ and
$\|\Gamma\|_{\Omega_{\Theta}}<\epsilon$, then
\begin{equation}\label{inequalitylocalstability}
\|N^{\overline{\mu},\Gamma}(t,.)\|\leq ce^{-d(t-\overline{t})},
\quad \textrm{ for } t\geq \overline{t}.
\end{equation}
\end{theorem}

\begin{proof}
Let $\varepsilon>0$ and $L\geq0$ be given by
(\ref{Lipshitzconditionneighborhoodzero}). We assume that
\begin{displaymath}
\|\overline{\mu}\|\leq \varepsilon \quad \textrm{ and } \quad \|
\Gamma \|_{\Omega_{\Theta}}\leq \varepsilon L.
\end{displaymath}
Then, Proposition \ref{invariance} implies that
\begin{displaymath}
|N^{\overline{\mu},\Gamma}(t,m)|\leq \varepsilon, \quad
\textrm{for } t\geq0 \textrm{ and } m\in[0,1].
\end{displaymath}
We define the sequence $(\overline{N}_n)_{n\in\mathbb{N}}$ by
\begin{displaymath}
\overline{N}_0(t,m)=K(t-\tau_{max},m)N^{\overline{\mu},\Gamma}\Big(\tau_{max},\chi\big(-(t-\tau_{max}),m\big)\Big),
\end{displaymath}
and
\begin{displaymath}
\overline{N}_n(t,m)=\overline{N}_0(t,m)-I_{\tau_{max}}(\overline{N}_{n-1})(t,m)+J_{\tau_{max}}(\overline{N}_{n-1})(t,m),
\end{displaymath}
for $t\geq\tau_{max}$, $m\in[0,1]$ and $n\in\mathbb{N}^*$, with
\begin{displaymath}
I_{\tau_{max}}(N)(t,m)=\displaystyle\int_{\tau_{max}}^t \!\!
K(t-s,m)\beta\Big(\chi(-(t-s),m),N\big(s,\chi(-(t-s),m)\big)\Big)N\big(s,\chi(-(t-s),m)\big)
ds,
\end{displaymath}
and
\begin{displaymath}
\begin{array}{rcl}
J_{\tau_{max}}(N)(t,m)&=&\displaystyle\int_{\tau_{max}}^t \!\! K(t-s,m) \zeta(\chi(-(t-s),m))N\big(s-\tau(\Delta(\chi(-(t-s),m))),\Delta(\chi(-(t-s),m))\big)\times\\
&&\beta\Big(\Delta(\chi(-(t-s),m)),N\big(s-\tau(\Delta(\chi(-(t-s),m))),\Delta(\chi(-(t-s),m))\big)\Big)
ds,
\end{array}
\end{displaymath}
and
\begin{displaymath}
\overline{N}_n(t,m)=N^{\overline{\mu},\Gamma}(t,m),\quad \textrm{
for } (t,m)\in[0,\tau_{max}]\times[0,1] \textrm{ and }
n\in\mathbb{N}.
\end{displaymath}
First, it is easy to check, by using the same reasonning as in the
proof of Proposition \ref{invariance}, that the sequence
$(\overline{N}_n)_{n\in\mathbb{N}}$ satisfies
\begin{displaymath}
|\overline{N}_n(t,m)| \leq \varepsilon, \qquad \textrm{ for }
t\geq0, \ m\in[0,1] \textrm{ and } n\in\mathbb{N}^*,
\end{displaymath}
since
\begin{displaymath}
|\overline{N}_0(t,m)| \leq \varepsilon, \qquad \textrm{ for }
t\geq0 \textrm{ and } m\in[0,1].
\end{displaymath}
Secondly, the assumption
\begin{displaymath}
L(1+2\kappa)<\widetilde{\delta},
\end{displaymath}
implies that
\begin{displaymath}
L(1+\|\zeta\|)<\widetilde{\delta}.
\end{displaymath}
Therefore, there exists $\rho\in(0,\widetilde{\delta})$ such that
\begin{displaymath}
L<\frac{\widetilde{\delta}-\rho}{1+\|\zeta\|e^{\rho\tau_{\max}}}<\frac{\widetilde{\delta}}{1+\|\zeta\|}.
\end{displaymath}
Let $t\geq\tau_{max}$ and $m\in[0,1]$. Then, we get
\begin{displaymath}
|\overline{N}_0(t,m)| \leq
e^{-\widetilde{\delta}(t-\tau_{max})}\varepsilon \leq
e^{-\rho(t-\tau_{max})}\varepsilon \leq\varepsilon.
\end{displaymath}
Since
\begin{displaymath}
|\overline{N}_1(t,m)-\overline{N}_0(t,m)|\leq
|I_{\tau_{max}}(\overline{N}_{0})(t,m)|+|J_{\tau_{max}}(\overline{N}_{0})(t,m)|,
\end{displaymath}
then, the estimates
\begin{displaymath}
|I_{\tau_{max}}(\overline{N}_{0})(t,m)|\leq  \varepsilon L
\int_{\tau_{max}}^t
e^{-\widetilde{\delta}(t-s)}e^{-\rho(s-\tau_{max})} ds,
\end{displaymath}
and
\begin{displaymath}
|J_{\tau_{max}}(\overline{N}_{0})(t,m)|\leq  \varepsilon L
\|\zeta\|\int_{\tau_{max}}^t
e^{-\widetilde{\delta}(t-s)}e^{-\rho(s-2\tau_{max})} ds,
\end{displaymath}
yield to
\begin{displaymath}
|\overline{N}_1(t,m)-\overline{N}_0(t,m)|\leq \varepsilon
L(1+\|\zeta\|e^{\rho\tau_{\max}})e^{-\widetilde{\delta}t}e^{\rho\tau_{\max}}\int_{\tau_{max}}^t
e^{(\widetilde{\delta}-\rho)s} ds.
\end{displaymath}
By remarking that
\begin{displaymath}
e^{-\widetilde{\delta}t}e^{\rho\tau_{\max}}\int_{\tau_{max}}^t
e^{(\widetilde{\delta}-\rho)s} ds \leq
\frac{1}{\widetilde{\delta}-\rho}e^{-\rho(t-\tau_{\max})},
\end{displaymath}
we obtain
\begin{displaymath}
|\overline{N}_1(t,m)-\overline{N}_0(t,m)|\leq \varepsilon
L\frac{1+\|\zeta\|e^{\rho\tau_{\max}}}{\widetilde{\delta}-\rho}e^{-\rho(t-\tau_{\max})}.
\end{displaymath}
It is easy to see by induction, that
\begin{displaymath}
|\overline{N}_n(t,m)-\overline{N}_{n-1}(t,m)|\leq \varepsilon
\bigg(L\frac{1+\|\zeta\|e^{\rho\tau_{\max}}}{\widetilde{\delta}-\rho}\bigg)^ne^{-\rho(t-\tau_{\max})},
\quad \textrm{ for } n\in\mathbb{N}^*.
\end{displaymath}
Hence, we get
\begin{displaymath}
\begin{array}{rcl}
|\overline{N}_n(t,m)|&\leq&\displaystyle\sum_{i=1}^{n} |\overline{N}_i(t,m)-\overline{N}_{i-1}(t,m)| + |\overline{N}_0(t,m)|,\\
    &\leq&\displaystyle\sum_{i=1}^{n}\varepsilon \bigg(L\displaystyle\frac{1+\|\zeta\|e^{\rho\tau_{\max}}}{\widetilde{\delta}-\rho}\bigg)^i e^{-\rho(t-\tau_{\max})} + \varepsilon e^{-\rho(t-\tau_{max})},\\
    &\leq&\varepsilon e^{-\rho(t-\tau_{\max})}\displaystyle\sum_{i=0}^{n}\bigg(L\displaystyle\frac{1+\|\zeta\|e^{\rho\tau_{\max}}}{\widetilde{\delta}-\rho}\bigg)^i.
\end{array}
\end{displaymath}
Since
\begin{displaymath}
L\frac{1+\|\zeta\|e^{\rho\tau_{\max}}}{\widetilde{\delta}-\rho}<1,
\end{displaymath}
then, we finally obtain that
\begin{displaymath}
|N^{\overline{\mu},\Gamma}(t,m)|\leq
\frac{(\widetilde{\delta}-\rho)\varepsilon}{\widetilde{\delta}-\rho
-L(1+\|\zeta\|e^{\rho\tau_{\max}})}e^{-\rho(t-\tau_{\max})}, \quad
\textrm{ for all } t\geq\tau_{max}.
\end{displaymath}
Hence, we have obtained the inequality
(\ref{inequalitylocalstability}) with $\overline{t}=\tau_{max}$,
$\epsilon=\max\{\varepsilon, \varepsilon L \}$, $d=\rho$ and
\begin{displaymath}
c=\frac{(\widetilde{\delta}-\rho)\varepsilon}{\widetilde{\delta}-\rho
-L(1+\|\zeta\|e^{\rho\tau_{\max}})}.
\end{displaymath}
This completes the proof.
\end{proof}

By using Theorem \ref{theoremlocalstabilityN}, we can prove the
local exponential stability of the trivial solution of Problem
(\ref{integratedformulationN})-(\ref{integratedformulationP}).
This result is presented in the following corollary.

\begin{corollary}\label{corlocalstab}
Assume that
\begin{displaymath}
L(1+2\kappa)<\min\{\widetilde{\gamma},\widetilde{\delta}\}.
\end{displaymath}
Then, the trivial solution of Problem
(\ref{integratedformulationN})-(\ref{integratedformulationP}) is
locally exponentially stable.
\end{corollary}

\begin{proof}
First, since $L(1+2\kappa)<\widetilde{\delta}$, then Theorem \ref{theoremlocalstabilityN} implies that the trivial solution of Equation (\ref{integratedformulationN}) is locally exponentially stable. Moreover, since we suppose that $L(1+2\kappa)<\widetilde{\gamma}$, we can show, by using the same arguments as in the proof of Theorem \ref{theoremlocalstabilityN} and a variation of constant formula for $P$, that the trivial solution of Equation (\ref{integratedformulationP}) is locally exponentially stable, as soon as the trivial solution of Equation (\ref{integratedformulationN}) has this property.
\end{proof}

Finally, by using Theorem  \ref{theoremlocalstabilityN}, we can
obtain the global stability of the trivial solution of Problem
(\ref{integratedformulationN})-(\ref{integratedformulationP}).
This is done in the next proposition.

\begin{proposition}
Assume that the mapping $x\mapsto x\beta(m,x)$ is Lipschitz
continuous for all $m\in[0,1]$, with a Lipschitz constant $L$, and
that
\begin{displaymath}
L(1+2\kappa)<\min\{\widetilde{\gamma},\widetilde{\delta}\}.
\end{displaymath}
Let $\Gamma\in C(\Omega_{\Theta})$ and $\mu$ be a continuous
function such that $\overline{\mu}\in C[0,1]$. If $\| \Gamma
\|_{\Omega_{\Theta}}\leq L\|\overline{\mu}\|$, then, the trivial
solution of Problem
(\ref{integratedformulationN})-(\ref{integratedformulationP}) is
globally exponentially stable.
\end{proposition}

\begin{proof}
We first show, as in the proof of Proposition \ref{invariance},
that
\begin{displaymath}
|N^{\overline{\mu},\Gamma}(t,m)|\leq\|\overline{\mu}\| , \quad
\textrm{for } t\geq0 \textrm{ and } m\in[0,1].
\end{displaymath}
Then, by using the same reasonning as in the proofs of Theorem
\ref{theoremlocalstabilityN} and Corollary \ref{corlocalstab}, we
conclude.
\end{proof}

We can give some explanations about the condition $L(1+2\kappa)<\min\{\widetilde{\gamma},\widetilde{\delta}\}$. In fact, we can notice that this inequality is satisfied if $\widetilde{\gamma}$ and $\widetilde{\delta}$ are large or if $L$ is small enough. This corresponds, biologically, to the case where the mortality rates ($\gamma$ and $\delta$) are important or to the case when only a few cells are introduced in the proliferating phase, and then the cells supply is not sufficient ($L$ is a bound of the number of introduced cells).

\section{Discussion}

It is usually believed that the function $\beta$ is a Hill
function (see Mackey \cite{mackey1978}) given by
\begin{displaymath} 
\beta(m,x)=\left\{ \begin{array}{ll}
\beta_{0}(m)\displaystyle\frac{\theta^{n}(m)}{\theta^{n}(m)+x^{n}},& \textrm{ for } x\geq 0, \ m\in[0,1],\\
\beta_{0}(m),& \textrm{ for } x< 0, \ m\in[0,1],
\end{array} \right.
\end{displaymath}
with $\theta$ and $\beta_0$ two continuous and positive functions,
and $n\geq 1$. In this case, the mapping $x\mapsto x\beta(m,x)$ is
always Lipshitz continuous, with a Lipschitz constant equals to
\begin{displaymath}
\sup_{m\in [0,1]} \beta_0(m).
\end{displaymath}
Hence, the results of Theorem \ref{theoremlocalexistence} and
Corollary \ref{globalsolution} hold.

By the same way, Condition (\ref{betapropertyfeedback}) is easily
satisfied, because $\beta(m,x)=\beta(m,0)$ if $x<0$, and, if
$x\geq 0$, the function $x\mapsto \beta(m,x)$ is decreasing for
all $m$. Hence, the positivity of the solutions is naturally
obtained.

As it has already been noticed by Dyson \emph{et al.}
\cite{webb1996}, Mackey and Rudnicki \cite{mackey1999} and Adimy
and Pujo-Menjouet \cite{adimypujo}, we can expect to show the
influence of the immature cells population (that means, the
population with a small maturity) over the entire population. This
has also been obtained by Adimy and Crauste \cite{adimycrauste}
for a model with a proliferating phase duration distributed
according to a density.

In particular, by using the result of local stability obtained in
Section \ref{sectionlocalstability}, we will certainly be able to
prove that the stability (or the instability) of the immature
cells population leads to the global stability (or the
instability) of the entire population. This has been displayed for
the first time by Mackey and Rudnicki \cite{mackey1999} in 1999.
Since the stem cells population is known to be at the root of the
blood production system, then, this behaviour is naturally
expected in our model, and it is the purpose of a next work.




\begin{thebibliography}{99}

\bibitem{adimycrauste} M. Adimy and F. Crauste, \emph{Un mod\`ele non-lin\'eaire de prolif\'eration cellulaire : extinction des cellules et invariance}, C. R. Acad. Paris Ser I \textbf{336}, 559-564 (2003).
%
\bibitem{adimycrauste2003} M. Adimy and F. Crauste, \emph{Global stability of a partial differential equation with distributed delay due to cellular replication}, Nonlinear Analysis {\bf 54}, 1469-1491 (2003).
%
\bibitem{adimypujo} M. Adimy and L. Pujo-Menjouet, \emph{A singular transport model describing cellular division}, C. R. Acad. Sci. Paris Ser. I Math. {\bf 332}, 12, 1071-1076 (2001).
%
\bibitem{adimypujo2} M. Adimy and L. Pujo-Menjouet, \emph{Asymptotic behaviour of a singular transport equation modelling cell division}, Dis. Cont. Dyn. Sys. Ser. B \textbf{3}, 3, 439-456 (2003).
%
\bibitem{adimypujo2003} M. Adimy and L. Pujo-Menjouet, \emph{A mathematical model describing cellular division with a proliferating phase duration depending on the maturity of cells}, accepted in Electron. J. Differ. Equ.
%
\bibitem{bradford} G. Bradford, B. Williams, R. Rossi and I. Bertoncello, \emph{Quiescence, cycling, and turnover in the primitive haematopoietic stem cell compartment}, Exper. Hematol. {\bf 25}, 445-453 (1997).
%
\bibitem{burns} F.J. Burns and I.F. Tannock, \emph{On the existence of a $G_{0}$ phase in the cell cycle}, Cell. Tissue Kinet. {\bf 19}, 321-334 (1970).
%
\bibitem{crabb1996_1} R. Crabb, J. Losson and M.C. Mackey, \emph{Dependence on initial conditions in non local PDE's and heredetary dynamical systems}, Proc. Inter. Conf. Nonlin. Anal. {\bf 4} (Tampa Bay, de Gruyter, Berlin), 3125-3136 (1996).
%
\bibitem{crabb1996_2} R. Crabb, M.C. Mackey and A. Rey, \emph{Propagating fronts, chaos and multistability in a cell replication model}, Chaos {\bf 6}, 477-492 (1996).
%
\bibitem{webb1996} J. Dyson, R. Villella-Bressan and G.F. Webb, \emph{A singular transport equation modelling a proliferating maturity structured cell population}, Can. Appl. Math. Quart. {\bf 4}, 65-95 (1996).
%
\bibitem{webb2000} J. Dyson, R. Villella-Bressan and G.F. Webb, \emph{A nonlinear age and maturity structured model of population dynamics. I : Basic theory.}, J. Math. Anal. Appl. {\bf 242}, 1, 93-104 (2000).
%
\bibitem{webb2000_2} J. Dyson, R. Villella-Bressan and G.F. Webb, \emph{A nonlinear age and maturity structured model of population dynamics. II : Chaos.}, J. Math. Anal. Appl. {\bf 242}, 2, 255-270 (2000).
%
%
\bibitem{john} P.C.L John, \emph{The cell cycle}, London, Cambridge University Press (1981).
%
\bibitem{kato2000} Kato N., \emph{Positive global solutions for a general model of size-dependent population dynamics}, Abstr. Appl. Anal. \textbf{5}, 3, 191-206 (2000).
%
\bibitem{mackey1978} M.C. Mackey, \emph{Unified hypothesis of the origin of aplastic anaemia and periodic hematopoiesis}, Blood {\bf 51}, 941-956 (1978).
%
%
\bibitem{mackey1993} M.C. Mackey and A. Rey,  \emph{Multistability and boundary layer development in a transport equation with retarded arguments}, Can. Appl. Math. Quart. {\bf 1}, 1-21 (1993).
%
\bibitem{mackey1995_1} M.C. Mackey and A. Rey, \emph{Transitions and kinematics of reaction-convection fronts in a cell population model}, Physica D {\bf 80}, 120-139 (1995).
%
\bibitem{mackey1995_2} M.C. Mackey and A. Rey, \emph{Propagation of population pulses and fronts in a cell replication problem : non-locality and dependence on the initial function}, Physica D {\bf 86}, 373-395 (1995).
%
\bibitem{mackey1994} M.C. Mackey and R. Rudnicki, \emph{Global stability in a delayed partial differential equation describing cellular replication}, J. Math. Biol. {\bf 33}, 89-109 (1994).
%
\bibitem{mackey1999} M.C. Mackey and R. Rudnicki, \emph{A new criterion for the global stability of simultaneous cell replication and maturation processes}, J. Math. Biol. {\bf 38}, 195-219 (1999).
%
\bibitem{mitchison} J.M Mitchison, \emph{The biology of the cell cycle}, London, Cambridge University Press (1971).
%
\bibitem{sachs} Sachs L., \emph{The molecular control of hemopoiesis and leukomia}, C. R. Acad. Sci. Paris \textbf{316}, 882-891 (1993).
%
\bibitem{traviswebb} Travis C.C. and Webb G.F., \emph{Existence and stability for partial functionnal differential equations}, Trans. Am. Math. Soc. \textbf{200}, 395-418 (1974).
%
\bibitem{webb1985} Webb G.F., \emph{Theory of non-linear age-dependent population dynamics}, Monographs and Textbook in Pure and Applied Mathematics, 89, New-York Basel: Marcel Dekker Inc., 294p, 1985.
%
%
\end{thebibliography}
\end{document}